\documentclass[10pt]{amsart}
\usepackage{calc}
\usepackage{graphicx}
\usepackage{tikz}
\usepackage{amssymb}
\usepackage{hyperref}
\usepackage[all]{xy}
\usepackage{colortbl}

\usepackage{color}
\definecolor{gray}{rgb}{.75,.75,.75}
\definecolor{orange}{rgb}{.9,.6,.2}
\definecolor{violet}{rgb}{.5,0,.5}
\definecolor{dg}{rgb}{0,0.67,0}
\definecolor{lg}{rgb}{0.7,0.7,0.7}
\definecolor{brown}{rgb}{0.58,0.29,0}

\newcommand{\nc}{\newcommand}
\nc{\bl}[1]{{\color{blue} #1}}
\nc{\br}[1]{{\color{brown} #1}}
\nc{\gr}[1]{{\color{dg} #1}}
\nc{\darkgr}[1]{{\color{darkgreen} #1}}
\nc{\re}[1]{{\color{red} #1}}
\nc{\ora}[1]{{\color{orange} #1}}
\nc{\gray}[1]{{\color{lg} #1}}
\nc{\ita}[1]{\color{red} {\mathit #1}}
\nc{\cR}{\mathcal{R}}

\theoremstyle{plain}
\newtheorem{theorem}{Theorem}[section]

\newtheorem{algorithm}[theorem]{Algorithm}

\theoremstyle{definition}

\DeclareMathOperator{\Aut}{Aut}
\DeclareMathOperator{\Amb}{Amb}
\DeclareMathOperator{\Gal}{Gal}

\newcommand{\Q}{\mathbb{Q}}

\newcommand{\R}{\mathbb{R}}

\newcommand{\C}{\mathbb{C}}
\newcommand{\Z}{\mathbb{Z}}
\newcommand{\Qp}{\mathbb{Q}_p}
\newcommand{\Qpsep}{\mathbb{Q}^{\rm sep}_p}

\newcommand{\F}{\mathbb{F}}
\newcommand{\Eisen}{\mbox{Eisen}}

\newcommand{\EisenII}{\mbox{Eis}}

\newcommand{\ord}{\mbox{ord}}

\newcommand{\cO}{\mathcal{O}}

\newcommand{\cmmt}[1]{}

\newcommand{\Inv}{\mbox{Inv}}

\numberwithin{figure}{section}
\numberwithin{equation}{section}
\numberwithin{table}{section}

\newcommand{\field}[1]{\href{https://lmfdb.org/padicField/#1}{\texttt{#1}}}
\newcommand{\ffield}[2]{\href{https://lmfdb.org/padicField/#1#2}{\texttt{#2}}}
\newcommand{\family}[1]{\href{https://lmfdb.org/padicField/family/#1}{\texttt{#1}}}
\newcommand{\ffamily}[2]{\href{https://lmfdb.org/padicField/family/#1#2}{\texttt{#2}}}

\newcommand{\kk}[1]{\marginpar{\parbox{3.5cm}{\scriptsize \color{blue} \sf KK:\  #1}}}

\newcommand{\jg}[1]{\marginpar{\parbox{3.5cm}{\scriptsize \color{blue} \sf JG:\  #1}}}

\title{Families of $p$-adic fields}
\author[J. Guardia]{Jordi Gu\`ardia-Rubies }
\address{Universitat Polit\`ecnica de Catalunya}
\email{jordi.guardia-rubies@upc.edu}

\author[J.W. Jones]{John W. Jones}
\address{Arizona State University}
\email{jj@asu.edu}

\author[K. Keating]{Kevin Keating}
\address{University of Florida}
\email{keating@ufl.edu}

\author[S. Pauli]{Sebastian Pauli}
\address{University of North Carolina Greensboro}
\email{s\_pauli@uncg.edu}

\author[D.P. Roberts]{David P.\ Roberts}
\address{University of Minnesota Morris}
\email{roberts@morris.umn.edu}

\author[D. Roe]{David Roe}
\address{Massachusetts Institute of Technology}
\email{roed@mit.edu}
\date{\today}

\begin{document}

\begin{abstract} We improve the database of $p$-adic fields 
in the LMFDB by systematically using Krasner-Monge 
polynomials and working relatively as well as absolutely.  These improvements organize 
$p$-adic fields into families.  They thereby make long lists of fields 
more manageable and various theoretical structures more evident.
In particular, the database now includes all degree $n$ extensions of $\Q_p$,
for $p<200$ and $n \leq 23$.
\end{abstract}

\maketitle

\setcounter{tocdepth}{1}

\vspace{-.1in}
\tableofcontents

\vspace{-.2in}

\section{Introduction}
\label{sec introduction}
In arithmetic parts of mathematics, it is often useful to work one prime at a time.  
When working at a single prime $p$, the field $\Q_p$ of $p$-adic numbers
commonly plays a central role.  Also important are finite-degree field 
extensions of $\Q_p$.  The number of isomorphism classes of such extensions 
of a given degree $n$ is finite and given by a formula due
to Monge \cite[Thm 1]{MoCount}.   Some cardinalities are given in Table~\ref{tab introtable}.
\begin{table}[htb]
\vspace{0 in}
\[
{\renewcommand{\arraycolsep}{3.5pt}
  \begin{array}{r|rrrrrrrrrrrrrrrr}
    n & 1 & 2 & 3 & 4 & 5 & 6 & 7 & 8 & 9 & 10 & 11 & 12 & 13 & 14 & 15 &
      16 \\
      \hline
   p= 2 & \gray{1} & 7 & \gray{2} & 59 & \gray{2} & 47 & \gray{2} & 1823 & \gray{3} & 158 & \gray{2} & 5493 & \gray{2} & 590 & \gray{4} &
      890111 \\
   p= 3 & \gray{1} & \gray{3} & 10 & \gray{5} & \gray{2} & 75 & \gray{2} & \gray{8} & 795 & \gray{6} & \gray{2} & 785 & \gray{2} & \gray{6} & 1172 & \gray{13}
      \\
    p=5 & \gray{1}& \gray{3} & \gray{2} & \gray{7} & 26 & \gray{ 7} & \gray{2} & \gray{11} & \gray{3}& 258 & \gray{2} & \gray{17} & \gray{2}& \gray{6} & 1012 & \gray{17 } \\
\hline
  p=2&   \gray{1} & 3& \gray{2} & 10 & \gray{2} & 8& \gray{2} & 49& \gray{3} & 10& \gray{2} & 43 & \gray{2} & 
 12 & \gray{4} & 389 \\
  p=3&   \gray{1} & \gray{2} & 4& \gray{3} & \gray{2} & 10& \gray{2} & \gray{4} & 28 & \gray{4} & \gray{2} & 20& \gray{2} & 
  \gray{4} & 16 & \gray{5} \\
  p=5 &   \gray{1} & \gray{2} & \gray{2} & \gray{3} & 6 & \gray{4} & \gray{2} & \gray{4} & \gray{3} & 16 & \gray{2} & \gray{6} & \gray{2} & 
  \gray{4} & 20 & \gray{5} 
    \end{array}
   }
 \]
 \caption{\label{tab introtable} Top: The number of isomorphism classes of degree
 $n$ field extensions of $\Q_p$.  Bottom: The number of degree $n$ families
 over $\Q_p$.  Cases where all ramification is tame are in gray.}
 \end{table}
 \vspace{-.2in}
 
 An online database giving defining polynomials and various invariants of $p$-adic fields 
 appeared in 2006, in connection with the paper \cite{jones-roberts-database}, which built on
 \cite{amano,panayi,pauli-roblot}.
 The original database quickly expanded via various works \cite{nonic,octic,awtrey-p3d12,awtrey-p2d12c4},
 to include complete detailed
 tables for many $(p,n)$, including all those listed on Table~\ref{tab introtable} except $(2,16)$.  
Another improvement was a migration in 2011 from an {\em ad hoc} platform
 to the LMFDB \cite{LMFDB}, so that the
 data can be more easily inspected from a wider variety
 of perspectives.  
 
The purpose of this paper is to describe a 
substantial improvement we have recently made to 
the database in the LMFDB.  The starting point for the
improvement is the systematic use of certain Eisenstein polynomials.  These polynomials 
were introduced long ago by Krasner \cite{krasner-cluj} 
and their theory was brought into modern form by Monge \cite{Mo}.
In this approach, fields are naturally organized 
 into families.   The number of families for small $(p,n)$ 
 is also given in Table~\ref{tab introtable}.  
 
 The improved database is at
 \begin{center}
 \url{https://lmfdb.org/padicField/}
 \end{center}
 and has a page for each of the 1,335,301 fields with $p < 200$ and $n \leq 23$. 
 The new case on Table~\ref{tab introtable}, $(p,n)=(2,16)$, by itself accounts for
 about $64\%$ of the fields in this range.  The next two largest contributors are
 $(2,20)$ and $(3,18)$, and they contain 314,543 and 130,647 fields respectively.
 These two
 new
  cases  account for about $23\%$ and $9\%$ of the total, so that 
 all the other $(p,n)$ together contribute about $4\%$.   
 
  The improved
 database also now has a page for each of the $19585$ families with
 $p<200$ and $n \leq 47$.
 It is not reasonable to populate all of these families with fields, 
 as there are in total about 897 billion.
 However, if a need arises for including all the fields belonging to a specific 
 reasonably-sized family,
 doing so will be straightforward.
 Some supporting code for the database is at \url{https://github.com/roed314/padic_db}.

 The subsections of this introduction give a first idea of 
 the previous structure of the database and how the systematic introduction
 of families improves it.  Sections~\ref{sec herbrand} and \ref{sec eisenstein} then present the theory necessary
 for the improvement, with some of the more subtle details and
 various algorithmic issues deferred to the companion paper \cite{polred}. 
Sections~\ref{sec sample} and \ref{sec theoretical} encourage the reader to explore and appreciate the database,
first by focusing on sample families and then by focusing on connections
to various theorems in the Galois theory of $p$-adic fields.

\subsection{The previous field-by-field approach} 
\label{subsec previous} We begin by describing 
some aspects of the database as it stood before our recent improvements.  
Degree $n$ fields were presented as $L=\Q_p[x]/f(x)$ with $f(x) \in \Z[x]$ 
a degree $n$ polynomial obtained from a search over possibilities.
Detailed attention was not paid to the choice of $f(x)$.  Rather the 
focus was on the most important invariants of $L$.

The extension $L/\Q_p$ has a normal closure 
$L^{\rm gal}/\Q_p$ and hence a Galois group
$\Gal(L^{\rm gal}/\Q_p)$.  The general theory of
$p$-adic fields gives a decreasing filtration of this 
group by normal subgroups.  The successive
subquotients $Q^s$ each have a size $|Q^s|$ and also 
an associated slope $s$.  Here slopes
of $-1$, $0$, and positive rational numbers 
correspond to no ramification, 
tame ramification, and wild ramification. 

The database focused on the filtered group
$\Gal(L^{\rm gal}/\Q_p)$.  Filtered groups
are somewhat unwieldy objects, so the
database gave only associated numerical
invariants.  To represent the group,
it gave its standard label $nTj$ in the list
of conjugacy classes of transitive
subgroups of the symmetric group 
$S_n$ \cite{butler-mckay,conway-hulpke-mckay,hulpke05}.  To represent the filtration, it
gave the {\em Galois slope content}
$W_t^u$.   Here $u=|Q^{-1}|$, $t=|Q^0|$,
and a wild subquotient $Q^s$ of size $p^\rho$ 
contributes $\rho$ copies of $s$ to the
weakly increasing list $W$ of wild 
slopes.  So the word ``content'' is
in the spirit of ``Jordan-H\"older content.''

As an example, one of the $795$ nonic $3$-adic fields $L$
was represented by the polynomial
\begin{equation}
\label{eq oldnonic}
f(x) = 21 + 18 x + 18 x^2 + 21 x^3 + 9 x^4 + 18 x^5 + x^9.
\end{equation}
The Galois group has $324 = 2^2 3^4$ elements.  On
he standard list from $9T1 = C_9$ to $9T45 = S_9$, it
is $9T24$.  The Galois slope content is $[\frac{1}{2},\frac{1}{2},\frac{2}{3},\frac{3}{2}]_2^2$. 
In general, Galois groups were determined in the above-cited papers by computing and factoring many resolvents over $\Q_p$.
The Galois slope content was determined by studying the ramification in the fields defined by these factors.   

Previous to our current improvements, the
database emphasized Artin slopes $\hat{s}_k = s_k+1$
rather than our current slopes $s_k$, displaying 
e.g.\ $[\frac{3}{2},\frac{3}{2},\frac{5}{3},\frac{5}{2}]_2^2$ in the above example.  
 Artin slopes $\hat{s}_k$ are indeed often
more convenient in global applications.  
We are now emphasizing the smaller
slopes $s_k$, often called Swan slopes, 
because they are
more natural in detailed local
analyses.

\subsection{The new family approach} 
\label{subsec new} 
The Galois slope content 
of any finite extension $L/\Q_p$ splits cleanly into the {\em visible slope content}
$[s_1,\dots,s_w]_\epsilon^f$, with associated degree
$[L:\Q_p] = f \epsilon p^w$, and the rest,  called {\em hidden slope content},
with associated degree $[L^{\rm gal}:L]$.  In  example \eqref{eq oldnonic}, 
the visible slope content is $[\frac{1}{2},\frac{3}{2}]_1^1$ and the 
hidden slope content is $[\frac{1}{2},\frac{2}{3}]_2^2$.   The visible slope
content is enormously easier to compute than the
hidden slope content.   As we will explain, 
no Galois-theoretic concepts are needed. 

We say that two extensions of $\Q_p$ belong
to the same {\em absolute family} if their visible slope
contents are the same.  The key idea underlying this paper
is that one can find defining polynomials 
for all the extensions in a family by suitably 
specializing a single {\em generic polynomial}
belonging to the family.  These specializations
are the previously mentioned Krasner-Monge
polynomials.  

Continuing the example begun in \eqref{eq oldnonic}, consider the family of $3$-adic fields 
with visible slope content $[\frac{1}{2},\frac{3}{2}]_1^1$.  Following the general 
recipe we will present, the generic polynomial is 
\begin{equation}
\label{eq nonicgenpoly}
f(a_3,a_{10},b_{11},b_{13},\pi,x) = \pi \left(1 + a_{10} \pi x  + b_{11}  \pi x^2 +  a_3  x^3  +  b_{13} \pi x^4  \right) + x^9.
\end{equation}
Specializing via $\pi=3$, $a_\sigma \in \{1,2\}$, and $b_\sigma \in \{0,1,2\}$ 
gives thirty-six 3-Eisenstein polynomials in $\Z[x]$.  
They represent bijectively the thirty-six entries on the 
list of $795$ nonic $3$-adic fields which have 
visible slope content $[\frac{1}{2},\frac{3}{2}]_1^1$.   It couldn't be 
easier!  Moreover, as an important bonus, the 
coordinates provided by generic polynomials 
often give rise to clean descriptions of 
the hidden invariants.  In the case 
of \eqref{eq nonicgenpoly} there are ten possibilities
for the pair consisting of the Galois group and
the hidden slope content, one pair 
being the above $(9T24,[\frac{1}{2},\frac{2}{3}]_2^2)$.  
Three fields have this pair, namely the
ones with $(a_3,a_{10},b_{11}) = (1,2,2)$.
The particular field defined by  \eqref{eq oldnonic} comes from 
$b_{13}=2$.   The other nine subsets likewise 
have very elementary descriptions, as can be seen 
in Table~\ref{tab nonicpackets}.

Our sample visible slope content $[\frac{1}{2},\frac{3}{2}]_1^1$ has two simplifying features: its unramified degree $f$ is $1$ and the
other data in the visible slope content measuring ramification is {\em rigid}, as we will explain in \S\ref{subsec rigidity}.  We broadly
describe the general case in this paper,
but defer a full treatment of the complications associated with
$f>1$ and nonrigidity to \cite{polred}.

Using symbols like $[\frac{1}{2},\frac{3}{2}]_1^1$ in a naming scheme for families would be unwieldy 
as part of a URL.
The database instead uses labels in the form $p.f.e.cL$, as in \family{3.1.9.18b} for the example.  Here $c$ is the common discriminant-exponent of all the fields in the family and the letter $L$ resolves ambiguity. Similarly, the mathematically ideal Eisenstein coefficients do not work well as labels identifying fields within a family.
The database instead appends a subfamily number
$\ell$ followed by a counter $j$, so that the example field \eqref{eq oldnonic} becomes
\field{3.1.9.18b2.4}.
In \eqref{eq nonicgenpoly}, two fields are in the same
subfamily if they have the same $a_{10}$.  The
general notion of subfamily involves residual polynomials
and is given in \cite{polred}.

\subsection{The new relative context}
\label{subsec relative}
The paper \cite{jones-roberts-database} and the subsequent papers 
extending the database were aggressively absolute: 
$p$-adic fields of degree $N$ were always given by a degree 
$N$ polynomial 
with coefficients in $\Z$.   However it is often
better to build fields in towers and the
theory of generic polynomials fits 
perfectly into this paradigm.  
The theory of slopes generalizes
to this relative context and we say
that two extensions $L_1/K$ and $L_2/K$  belong
to the same family over $K$ if they
have the same relative visible slope content 
$I = [s_1,\dots,s_w]_\epsilon^f$.    
We denote this family, viewed simply as a finite
set of isomorphism classes of 
extensions, by $I/K$. 

Henceforth, we call $I$ a 
{\em Herbrand invariant}, 
because we have other ways of expressing the data in $I$ that do not directly involve slopes, 
as we will be explaining in the next section.  In fact our viewpoint is that $p$-adic Herbrand invariants are simple combinatorial objects that could be described independently of $p$-adic fields.  To get a family, 
one combines two objects of different nature, $I$ and $K$, subject to a numerical compatibility 
condition.

To continue the example begun in \eqref{eq oldnonic} yet further, take any finite
extension $K$ of $\Q_3$ as 
ground field, with residual cardinality
denoted by $q$.    Consider all 
extensions $L/K$ with Herbrand invariant
$I = [\frac{1}{2},\frac{3}{2}]_1^1$.  
This family $I/K$ is again bijectively indexed by certain specializations
of the exact same generic polynomial \eqref{eq nonicgenpoly}.
The difference is that $\pi$ now is specialized to a uniformizer
of $K$, rather than to the uniformizer $3$ of $\Q_3$, and the $a_i$ and $b_i$ to elements of $K$ with distinct reductions modulo $\pi$.
Thus the real purpose of \eqref{eq nonicgenpoly} is to get 
all relative extensions $L/K$ of type 
$[\frac{1}{2},\frac{3}{2}]_1^1$, for any 
fixed $3$-adic base field $K$.   
Directly generalizing the case $K=\Q_3$,
this family has cardinality $|I/K| =(q-1)^2 q^2$.  

While we work relatively throughout this paper, 
the improved database 
 keeps the original 
context of $p$-adic fields $L$ as one
of its two focal points.    It has a basic 
bipartite structure.  On the one hand,
each $p$-adic field $L$ 
within range has, as before, a homepage.  
On the other hand, each family $I/K$ now also has
a homepage.  The two homepages are linked if $K$ 
can be realized as a subfield of $L$ such
that $L$ is in the family $I/K$.

In the case of the example family $[\frac{1}{2},\frac{3}{2}]_1^1/\Q_3$, each field $L$ in it has exactly one cubic subfield $K'$ having discriminant-exponent $3$. The
cubic extension $L/K'$ then appears in the relative family
$[\frac{7}{2}]_1^1/K'$.  There are two possibilities
for $K'$, namely \field{3.1.3.3a1.1} or \field{3.1.3.3a1.2},
each of which occurs for half of the $L$.  
So the fields in the absolute family $[\frac{1}{2},\frac{3}{2}]_1^1/\Q_3$ come half each from the relative families
\family{3.1.3.3a1.1-1.3.9a} and 
\family{3.1.3.3a1.2-1.3.9a}.  Here, the
syntax for relative families is (base field)-$f.e.c\mbox{(tiebreaker)}$. 
Two fields $L$ in $[\frac{1}{2},\frac{3}{2}]_1^1/\Q_3$ also have three more cubic subfields, now with
discriminant-exponent $5$.  The 
index page \href{https://lmfdb.org/padicField/families/?label_absolute=3.1.9.18b}{\texttt{lmfdb.org/padicField/families/?label\_absolute=3.1.9.18b}} gives 
an overview of all possibilities.

\subsection{Notation} 
\label{subsec notation} We gather and comment on our most basic notations for the reader's convenience.  As already indicated, a prime $p$ is
fixed, the symbol $K$ is reserved for a finite extension of the field $\Q_p$ of 
$p$-adic numbers, and the symbol $L$ is reserved for a
finite extension of $K$.  We let $\cO$ be the ring of integers
of $K$, $\Pi$ its maximal ideal, $\kappa = \cO/\Pi$ its residue
field, and $q = |\kappa|$ its residual cardinality.  
 We usually view $K$ as fixed and $L$ as varying.  

Many fields and numbers are associated to a given $L/K$.  The most basic come from its standard
tower
\begin{equation}
\label{eq standard}
K \stackrel{f} \subseteq L_{\rm ur} \stackrel{\epsilon}{\subseteq} L_0 \stackrel{p^w}{\subseteq} L.
\end{equation}
Here $L_{\rm ur}$ is the maximal unramified subextension and $L_0$ is the maximal tamely
ramified subextension.
 The superscripts indicate relative degrees so that the entire degree $n=[L:K]$ comes with a canonical
factorization, $f \epsilon p^w$. 
In many circumstances it is enough to work with  the ramification index $e = \epsilon p^w$,
but in this paper it is usually best to separately emphasize its tame factor $\epsilon$ and its wild factor $p^w$.

The most familiar quantity capturing ramification in $L/K$ is the discriminant 
$\mbox{disc}(L/K) = \Pi^c$.  We are emphasizing the importance of the
discriminant-exponent $c$ in our LMFDB labeling scheme.   But for deeper analysis we prefer to 
switch to the {\em mean} $m$ via the transformation equation
\begin{equation}
\label{eq meandef} c = f(e-1+em). 
\end{equation}
So in the continuing nonic $3$-adic example, the focus on $c=18$ is shifted onto
$m=10/9$.

Most of our attention is focused on a discrete invariant $W$ measuring the 
wild ramification present in the extension $L/L_{\rm ur}$.  
One way to describe $W$ is by the vector $[s_1,\dots,s_w]$ already
emphasized.    Two similar ways are given in \eqref{eq herbrandmeans} and
\eqref{eq herbrandrams} and relations are summarized either by
\eqref{eq slopetransform} and \eqref{eq ramtransform} or, in a different manner, by 
\eqref{eq scaled} and \eqref{eq twotransforms}.
The invariant $W$ is the wild part of 
the Herbrand invariant $I = \Inv(L/K) = W_\epsilon^f$.

The set of all extensions $L/K$ sharing a given $I$ is denoted $I/K$.  
We think of $I$ as a discrete invariant of $L/K$. 
In contrast, we view the invariants necessary to distinguish fields 
inside of $I/K$ as continuous.  They have their own complicated notation,
introduced in Section~\ref{sec eisenstein}.

All the notation for $L/K$ has its analog for $K/\Q_p$, starting with 
$\Q_p \subseteq K_{\rm ur} \subseteq K_0 
\subseteq K$. 
If we were to name everything, we would have to take care to 
avoid notation clashes.  
Fortunately, our considerations here make very little use of the internal
structure of $K$, as we are focused mainly on extensions of 
$K$.  The invariants of $K$ which enter our considerations most often are
its absolute ramification index $e_K = [K:K_{\rm ur}] = \ord_\Pi(p)$ and
the above-mentioned residual cardinality $q = |\kappa|$.

In Sections~\ref{sec herbrand} and \ref{sec eisenstein} we work completely 
constructively, making no mention of any algebraic closure of $\Q_p$.  
In Sections~\ref{sec sample} and \ref{sec theoretical} the attention
shifts to Galois theory.  There we fix a separable, i.e.\ algebraic,  closure $\Qpsep$
of $\Q_p$ and view the ground field $K$ as in $\Qpsep$.  Given 
an abstract extension $L/K$ of degree $n$, we let $L^{\rm gal}$ be the 
compositum of the images of the $n$ different $K$-linear embeddings
of $L$ into $\Qpsep$.  So the Galois group $\Gal(L^{\rm gal}/K)$
is a quotient of $\Gal(\Qpsep/K)$.

\subsection{The Laurent series case} There is an extremely parallel case where one replaces the 
absolute ground field $\Q_p$ by the field $\F_p((t))$ of Laurent series over 
$\F_p$.   Here a separable closure is smaller than an algebraic closure.  It is the former which is 
the setting for the analog, which accounts for our notation $\Qpsep$ in this paper. 
The direct connections between the two cases are of great current
interest \cite[\S2]{Scholze}.  
As the database contains only extensions of 
$\Q_p$, we will limit our discussion of the Laurent series case to occasional 
brief remarks that clarify the $p$-adic case.

\section{Herbrand invariants}
\label{sec herbrand}
   Herbrand invariants are fundamental to this paper because they index families.  This section explains what they are and how to work with them explicitly.

\subsection{An elementary approach emphasizing canonical subfields} Serre's standard text  \cite[IV]{Se} 
associates a Herbrand function $\phi: \R_{\geq 0} \rightarrow \R_{\geq 0}$
 to any Galois extension $L/K$ of $p$-adic fields.  
The construction involves filtering $\Gal(L/K)$ by a descending
family of subgroups indexed in two ways, by lower and upper 
numbering.    Remark 1 of \cite[IV.3]{Se} extends the definition of 
$\phi$ to general extensions $L/K$ by reduction to the Galois case.    
 Details of this reduction were first explained in print by Deligne in a 
 six-page appendix to \cite{De}.   
 
\begin{figure}[htb]
\begin{center}
\includegraphics[width=4in]{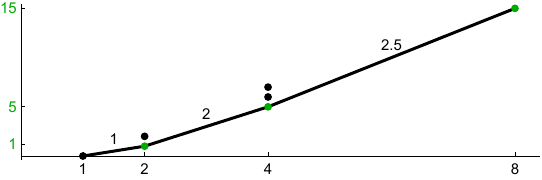}
\end{center}
\caption{\label{fig slopepolygon} The slope polygon $S$ associated to the
octic extension $L/\Q_2$ defined by \eqref{eq octic}.  The slopes
$[1,2,2.5]$
and \gr{means $\gr{\langle \frac{1}{2}, \frac{5}{4}, \frac{15}{8} \rangle = \langle 0.5,1.25,1.875 \rangle}$} 
are indicated. } 
\end{figure}

 It is less widely known that the Herbrand function
 of a general extension $L/K$ can be understood without any reference to group theory
 whatsoever.   In this subsection, we use the method of {\em slope polygons} to get 
 relevant numerical quantities.  We are 
 following  \cite[\S3.4]{jones-roberts-database}, except that the Artin slopes
 there are translated to Swan slopes here.    Figure~\ref{fig slopepolygon} starts our 
 second continuing example, based on 
  \begin{equation}
 \label{eq octic} f(x) = 1 + 6 x^4 + x^8.
 \end{equation}
 This particular polynomial is chosen because it defines a Galois 
 extension over both $\Q$ and $\Q_2$, with Galois group $D_4$; note that $f(x+1)$ is Eisenstein.
 The reader can then follow along, using number field software
 to confirm our statements, e.g.\ {\em Pari}'s \verb@nfsubfields@ 
 to get subfields and \verb@smalldiscf@ to get their discriminants
 and ultimately their means via \eqref{eq meandef}.   The fact that this
 example is Galois is irrelevant for the method we are describing.  

 Consider towers $L/L'/L_{\rm un}/K$.   
 For each, one has the degree $e' = [L':L_{\rm un}]$, the mean 
 $m'$ of $L'/L_{\rm un}$, and thus a point $(e',e'm')$ in a Cartesian plane.
 The {\em slope polygon} $S$ is the lower boundary of the
 convex hull of all such points.  Over each interval 
 $[\epsilon p^{k-1},\epsilon p^k]$ the polygon $S$ is just 
 a segment with some slope $s_k$.  The symbol 
 $[s_1,\dots,s_w]^f_\epsilon$ is then the visible slope 
 content emphasized in the introduction. 

 Define $h$ to be the function on $[1,e]$ having 
 graph $S$.    For $k=1$, \dots, $w$, define
 $m_k = \frac{h(\epsilon p^k)}{\epsilon p^k}$.  Clearly, the data 
 $\langle m_1,\dots,m_w \rangle^f_\epsilon$
 contains the exact same information as $[s_1,\dots,s_w]^f_\epsilon$.   
 If the point $P_k  = (\epsilon p^k,\epsilon p^k m_k)$ is
 a turning point or the right endpoint of $S$ then 
 we say that the index $k$ is final.  Then $P_k$
 comes from exactly one tower  
 $L/L_k/L_{\rm un}/K$.   The standard
 chain \eqref{eq standard} can be extended
 to a more refined chain from $K$ to 
 $L=L_w$ by including all the other 
 {\em canonical subfields} $L_k$.  Note that for 
 a final index $k$, the extension $L_k/L_{\rm un}$
 has mean $m_k$.  For a non-final index $k$, there
 may or may not be a extension mapping to $P_k$.  
 In the example, the indices $k=1$, $2$, and $3$ are all final.

 \subsection{Four perspectives on Herbrand invariants}
 The {\em slopes} $s_k$ just introduced are usually 
 called {\em breaks in the upper numbering}.  
 There is similarly a dual polygon called
 the ramification polygon in which 
 certain elementary quantities $r_k$ appear
 along with the means $m_k$ again.   These {\em rams} $r_k$
 coincide with the {\em breaks in the lower numbering}.  
 We will emphasize our elementary and
 non-Galois-theoretical viewpoint by using these
two terms systematically,
 and not making further reference to breaks 
 or numbering.  Even
 when we bring back Galois groups in the
 last two sections, we will use the briefer
 terminology of slopes.  
 
\label{subsec multiple}  We write the Herbrand invariant 
of a $p$-adic extension $L/K$ in four ways, the first 
classical but the others advantageous in various 
situations:
 \begin{align}
 \displaybreak[0]
\label{eq herbrand} \Inv(L/K) & =  (\phi,f) && \mbox{($\phi$ is the Herbrand function)} \\
\label{eq herbrandslopes}                & = [s_1,\dots,s_w]_\epsilon^f && \mbox{(the $s_k$ are the slopes)} \\
\label{eq herbrandmeans}             & = \langle m_1,\dots,m_w \rangle_\epsilon^f && \mbox{(the $m_k$ are the means)} \\
\label{eq herbrandrams}               & = (r_1,\dots,r_w)_\epsilon^f && \mbox{(the $r_k$ are the rams).}
 \end{align}
 As a matter of notation, the subscript $\epsilon$ and superscript $f$ are allowed to be omitted when they are $1$.   
 Our key reference \cite{Mo} had different aims that did 
 not require emphasis on Herbrand invariants.
 However it makes essential use of all the quantities 
 $s_k$, $m_k$ and $r_k$.  

\begin{figure}[htb]
 \begin{center}
 \includegraphics[width=4.5in]{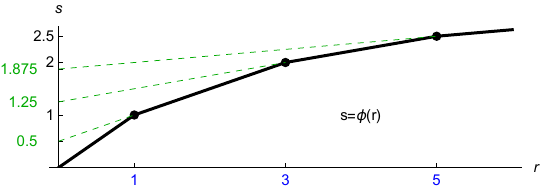}
 \end{center}
 \caption{\label{fig herbrandphi} The Herbrand function $\phi$
 for the octic extension $L/\Q_2$ defined by \eqref{eq octic}.  
 It takes the rams $\bl{(r_1,r_2,r_3) = (1,3,5)}$
 of this extension to its slopes 
 ${[s_1,s_2,s_3]=[1,2,2.5]}$.  The means $\gr{\langle 0.5,1.25,1.875 \rangle}$ are
 obtained by the indicated extensions of segments.} 
 \end{figure}

 As a convention, we put $s_0=m_0=r_0=0$.   Then \eqref{eq herbrandslopes}-\eqref{eq herbrandrams} are related via 
 \begin{align}
\label{eq slopetransform} s_k & =  {\frac{p m_k-m_{k-1}}{p-1}}, &   {m_k} & =  {\sum_{j=1}^k \frac{p-1}{p^{k+1-j}}  s_j}, \\
\label{eq ramtransform} {r_k} & = \epsilon p^k \frac{m_k- m_{k-1}}{p-1}, &  {m_k} & = {\sum_{j=1}^k \frac{p-1}{\epsilon p^j}  r_j}.  
\end{align}
Formulas \eqref{eq slopetransform} reflect the 
perspective of slope polygons:
each slope $s_k$ is a certain rise-over-run and each mean $m_k$ is a weighted average of 
slopes, with the formal slope $s_0=0$ having coefficient $1/p^w$, so that the coefficients sum to one.
Formulas \eqref{eq ramtransform} reflect the dual perspective of ramification
polygons.  Direct transformation formulas between slopes and rams
are a little more complicated, and one can just compose two transformations
with means in the middle.   

Figure~\ref{fig herbrandphi} uses our octic $2$-adic continuing example to illustrate how
the classical version \eqref{eq herbrand} is related to the more numerical versions
\eqref{eq herbrandslopes}-\eqref{eq herbrandrams}.
In general, the graph of $\phi$ starts at $(r_0,s_0) = (0,0)$ and 
goes linearly to the $(r_k,s_k)$ in order, 
with an actual step being taken only if 
$k$ is a final index, as otherwise  
we have $(a_{k},b_{k}) = (a_{k+1},b_{k+1})$.
It ends with a ray emanating from 
$(r_w,s_w)$.  If $k$ is $0$ or a final index then the slope of the segment with
left endpoint $(r_k,s_k)$ is $1/(\epsilon p^k)$.  With this strong condition on
slopes of the Herbrand segments, just $\phi$ determines 
all the $s_k$, $m_k$, and $t_k$.

\subsection{Automorphisms, mass, and rigidity}
\label{subsec rigidity} An automorphism in $\Aut(L/K)$ 
necessarily stabilizes all of the canonical subextensions
of $L/K$. Because of this fact, the Herbrand invariant $I$ alone
constrains the size of $\Aut(L/K)$.  A key input is that 
if $k$ is a final index, with the ram $r_k$ having been repeated $\rho$ times,
then the inertia group 
associated to $L_k/L_{k-\rho}$ has the form
$C_p^\rho \rtimes C_d$, where $d$ is the denominator
of the ram $r_k$.  The action is such that the step can only have nontrivial automorphisms if $d=1$, i.e.\ if $r_k$ is integral.

Define the {\em ambiguity number} of a $p$-adic Herbrand 
invariant $(r_1,\dots,r_w)_\epsilon^f$ to be 
$\Amb(I) = f \epsilon p^i$, 
where $i$ is the number of integral rams (see \S \ref{subsec reduction} for the definition of
$\Amb(I/K)$, the ambiguity of $I$ over $K$).
Then an extension $L/K$ in any $I/K$ has $|\Aut(L/K)|$ dividing $\Amb(I)$.  
The {\em mass} of $L/K$ is by definition $1/|\Aut(L/K)|$ and
$L/K$ is called {\em rigid} if its mass is $1$.  We say 
that $I$ is {\em rigid} if $\Amb(I)=1$.  So all 
extensions $L/K$ in a family $I/K$ with a rigid $I$ 
are rigid.

\subsection{Classification of Herbrand invariants} 
\label{subsec classification}
To classify all Herbrand invariants and say which actually occur over a given 
ground field $K$, it is best to use rams as was done in different language in 
\cite[Prop 3.10]{pauli-sinclair}.  
First consider totally wildly ramified extensions $L/K$ of degree $p^\rho$ having 
just a single ram $r$ repeated $\rho$ times.   The possibilities for $r$ depend 
only on the absolute ramification index $e_K$ of the ground field $K$.  

Define first $\cR_\rho^\infty$ to be the set of positive rational numbers with
denominator dividing $p^\rho-1$ and numerator not divisible by $p$.
Then $\cR_\rho^\infty$ is the correct collection of $r$ in the parallel case of 
Laurent series ground fields $\F_q((t))$.  
For the characteristic zero fields $K$ on which this paper is focused, the set of
possibilities $\cR_\rho^{e_K}$ contains  all the elements of $\cR_{\rho}^\infty$ which are less than $pe_K/(p-1)$.  These elements exhaust $\cR_\rho^{e_K}$ 
except that $\cR^{e_K}_1$ also contains $pe_K/(p-1)$.  We call this last ram 
{\em arithmetic} and all the smaller rams {\em geometric}.

Now consider all extensions $L/K$ with canonical factorization $f \epsilon p^w$. 
The Herbrand invariants that can arise are $(r_1,\dots,r_w)_\epsilon^f$  
where the $r_k$ are weakly increasing and each in their allowed set, from the
above considerations.   For example, to get all strictly increasing
sequences, each $r_k$ is chosen from $\cR_1^{e_K \epsilon p^{k-1}}$.     This description lets one
produce any entry on the lower half of Table~\ref{tab introtable}.  If $p$ is relatively prime to $n$, then the wild
considerations are all vacuous and the possible Herbrand invariants are simply $(\;)_\epsilon^{n/\epsilon}$ as
$\epsilon$ runs over divisors of $n$.

\begin{figure}[htb]
\begin{center}
\includegraphics[width=4.5in]{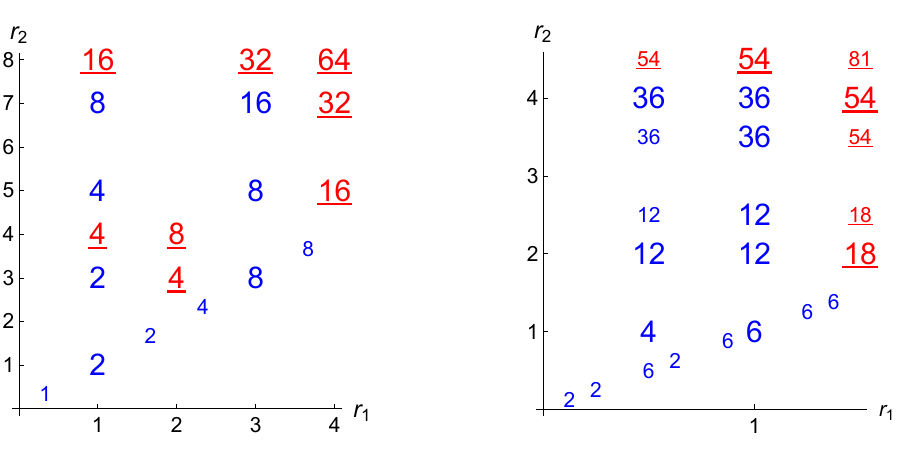}
\end{center}
\caption{\label{fig families}  $p$-adic Herbrand invariants $I=(r_1,r_2)$ with $r_j \leq k p^j/(p-1)$.  
 Left: $p = 2$ and $k = 2$.  Right: $p = 3$ and $k = 1$} 
\end{figure}

Figure~\ref{fig families} illustrates the case of Herbrand invariants $(r_1,r_2)$ over $2$- and $3$-adic fields,
with an eye towards giving a visual understanding of the general case.  Rather than represent a Herbrand invariant $I$ by a point at $(r_1,r_2)$, we
represent $I$ by an integer.  This integer is the mass of the family $I/K$, where $K$ is 
any ground field compatible with $I$ and having residue field of size $p$, as explained around \eqref{eq massformula}.
Families for which both rams are
geometric are presented in blue, 
and families which have an arithmetic
ram are presented in underlined red.  
Making an independent distinction, rigid
Hurwitz invariants are indicated by a smaller font. 

One should imagine each half of Figure~\ref{fig families} extended to the 
entire first quadrant of its $(r_1,r_2)$ plane, so that for each positive integer $k$ there is an underlined red 
hook $H_k$ of entries with upper right corner at $(k p/(p-1),k p^2/(p-1))$. 
 The Herbrand invariants $(r_1,r_2)$ compatible with a
 given $p$-adic field $K$ are exactly the arithmetic ones on $H_{e_K}$ and the
 geometric ones under $H_{e_K}$.  In other words, they are the ones on $H_{e_K}$, which are all underlined
 red, and the blue ones beneath $H_{e_K}$.
 
 As a numerical example, consider the total number of families over $\Q_p$ with degree 
 $p^2$.  There is just one family with ramification index $e=1$, the unramified family $(\;)^{p^2}$.  There
 are $p$ families with $e=p$, the geometric families $(1)^p$, $(2)^p$, \dots, $(p-1)^p$ and the
 arithmetic family $(p)^p$.  Guided by Figure~\ref{fig families}, one can check that there are
 $p^3 - \frac{p^2}{2} + \frac{p}{2}$ families with $e=p^2$.   For $p=2$ and $3$ respectively, 
 the total counts are $1+2+7=10$ and $1+3+24=28$, with $7$ and $24$ being the numbers
 of families on or under $H_1$ on Figure~\ref{fig families} and $10$ and $28$ appearing
 as entries in the lower half of Table~\ref{tab introtable}.  

 To see the seven and twenty-four families $p.1.p^2.C$ 
  listed as an indexing table, one can search for absolute families with residual characteristic $p$, residual field degree $f=1$, and ramification degree $e=p^2$. 
  To see
  a sixteen-line table corresponding to the $e_K=2$
  case of the left side of Figure~\ref{fig families}, one
  can search analogously for relative families.  While there are nineteen numbers in the 
  left side of the figure, the three in the underlined red hook do not correspond to families, 
  because they would involve an arithmetic ram, prohibited 
  by the above rules.  One needs to input an appropriate ground field by its label, say $K=\Q_2(i)$ by \field{2.1.2.2a1.1}, 
  as well as $f=1$ and $e=4$ again. 
  The three tables produced
  in this paragraph contain summary information, such as the masses appearing in Figure~\ref{fig families}.

\subsection{Explicit composition via sorting}
\label{subsec composition}
   Let $L/K/Q/\Qp$ be a tower with
 \begin{align}
 \label{eq twostep}
 \Inv(K/Q) & = (t_1,\dots,t_{w'})_{\epsilon'}^{f'} &
 \mbox{ and } && \Inv(L/K) & = (t_{w'+1},\dots,t_{w'+w''})_{\epsilon''}^{f''}.  
  \end{align}
 A great virtue of the traditional functional presentation of Herbrand invariants is that
 one has the simple formula $\phi_{L/Q}(r) = \phi_{K/Q}(\phi_{L/K}(r))$.
 Here we translate
 this formula into the more explicit language of rams, so as to 
 be able to go directly from \eqref{eq twostep} to $\Inv(L/Q) = (r_1,\dots,r_w)_\epsilon^f$
 
   Of course, $f=f' f''$, $\epsilon = \epsilon' \epsilon''$, and 
 $w = w'+w''$.  
To do the nontrivial part, we first  define $T_k = \epsilon'' t_k$ and formally write
 \begin{equation}
 \label{eq Tt}
 \Inv(L/Q) = (T_1,\dots,T_{w''},t_{w''+1},\dots,t_{w})_\epsilon^f.
 \end{equation}
 Then we change the entries of the $w$-vector by iterating the following replacement in
 any order.
Whenever there are two adjacent entries $(a,b)$ with $a>b$, replace them with $(b,b+p(a-b))$.  
 When the $w$-vector becomes weakly increasing, it is the desired
 ram vector $(r_1,\dots,r_w)$.  
 
The database uses this procedure to pass from the Herbrand invariant of a general relative
extension $L/K$
to the Herbrand invariant of the corresponding absolute
extension $L/\Q_p$. 
If $K$ is a canonical subfield of $L$, both steps of the process
are trivial and, very simply, $r_k = t_k$ for all $k$.  

To see the duality between rams and slopes, one can consider 
a modification of this process.   Step 1 yielding \eqref{eq Tt} is exactly 
the same.   In the modified Step 2, one iteratively replaces
adjacent increasing $(a,b)$ by $(b,b-(b-a)/p)$ until
the $w$-vector becomes weakly decreasing.  
As a new final Step 3, one reverses
the $w$-vector obtained to make it weakly increasing.  The
$k^{\rm th}$ entry of this final vector is $\epsilon s_k$, 
where $s_k$ is the $k^{\rm th}$ slope of the extension $L/Q$.

\section{Eisenstein polynomials} 
\label{sec eisenstein}
This section pictorially describes a well-behaved nonempty finite subset $\EisenII(L/K)$ of the infinite set $\Eisen(L/K)$ of all Eisenstein polynomials defining %CHANGE
any given totally ramified extension $L/K$.  Thus $f=1$ for the entirety of this section.

\subsection{Set-up}
\label{subsec setup} Let $\Eisen(e/K)$ be the 
space of Eisenstein polynomials of degree $e$ over $K$.  We write an element of $\Eisen(e/K)$ as  
\begin{equation}
\label{eq eisensteinF}
f(x) = F_0 + F_1 x + \dots + F_{e-1} x^{e-1} + x^e.
\end{equation}
So the $F_i$ run over the maximal ideal $\Pi$ of the ring of integers $\cO$ of $K$,
except that $F_0$ is not in $\Pi^2$.
One has decompositions into finitely many parts,
\begin{align}
\label{eq Eisenstart}
\Eisen(e/K) & =  \sqcup_I \Eisen(I/K),  & \Eisen(I/K) & = \sqcup_L \Eisen(L/K).
\end{align}
On the left, $I$ runs over totally ramified degree $e$ 
Herbrand invariants that are compatible with $K$.  On the
right, $L$ runs over isomorphism classes of extensions of $K$ which
have Herbrand invariant $I$.  

To present things as concretely as possible, we choose a generator $\pi$ of $\Pi$.  
We choose also a set of representatives $\tilde{\kappa}$ containing $0$ and $1$ for 
the $q$-element residue field $\kappa$.  
Rather than work with the $F_i$, we will work with $\pi$-adic 
expansions, writing each $F_i$ as  
$\pi \sum_{j=0}^\infty f_{i,j} \pi^j$, with $f_{i,j} \in \tilde{\kappa}$. 
If $K$ is unramified over $\Q_p$ we always take
the uniformizer $\pi$ to just be $p$.  A theoretically natural 
choice would be to take the nonzero elements of 
$\tilde{\kappa}$ to be the $(q-1)^{\rm st}$ roots of 
unity in $K$.  However, we make the computationally more convenient choice of $\tilde{\kappa} = \{0, \dots, q-1\}$ when $q$ is prime.

It will also be convenient to use single-indexing simultaneously with double-indexing, 
with $\sigma$ and $(i,j)$ determining each other
via $\sigma = j e+i$, $j = \lfloor \sigma/e \rfloor$, and $i \in \{0,\dots,e-1\}$.    With this convention, \eqref{eq eisensteinF}
is written
\begin{equation}
\label{eq eisensteinf}
F(x) = \pi \left( \sum_{\sigma=0}^\infty f_\sigma \pi^j x^i \right) + x^e.
\end{equation}
The last term $x^e$ will not be mentioned much in
our narrative because it does not contain an unspecified coefficient.  

A problem with the right part of \eqref{eq Eisenstart} is that the sets involved are infinite.  We will be replacing these sets by finite nonempty sets $\EisenII(I/K) = \sqcup_L \EisenII(L/K)$.  In the case when $I$ is rigid, each
of the $\EisenII(L/K)$ contains just one polynomial.

\subsection{Eisenstein diagrams} Fix a totally ramified Herbrand invariant
\label{subsec diagrams}
$I = [s_1,\dots,s_w]_\epsilon$.    To explicitly describe the 
sets $\EisenII(I/K)$ for all compatible $K$ at once, 
we consider the {\em Eisenstein diagram of $I$} in the vertical strip $R$ with 
horizontal coordinate $i \in [0,e)$ and vertical
coordinate $s \in [0,\infty)$.   As indicated by the name, this diagram depends 
on $I$ only, not on $K$.   Always we draw only $ [0,e) \times [0,s_w]$ 
as there is no useful information associated to the rest of the strip.  

\begin{figure}[htb]
\begin{center}
\includegraphics[width=4.5in]{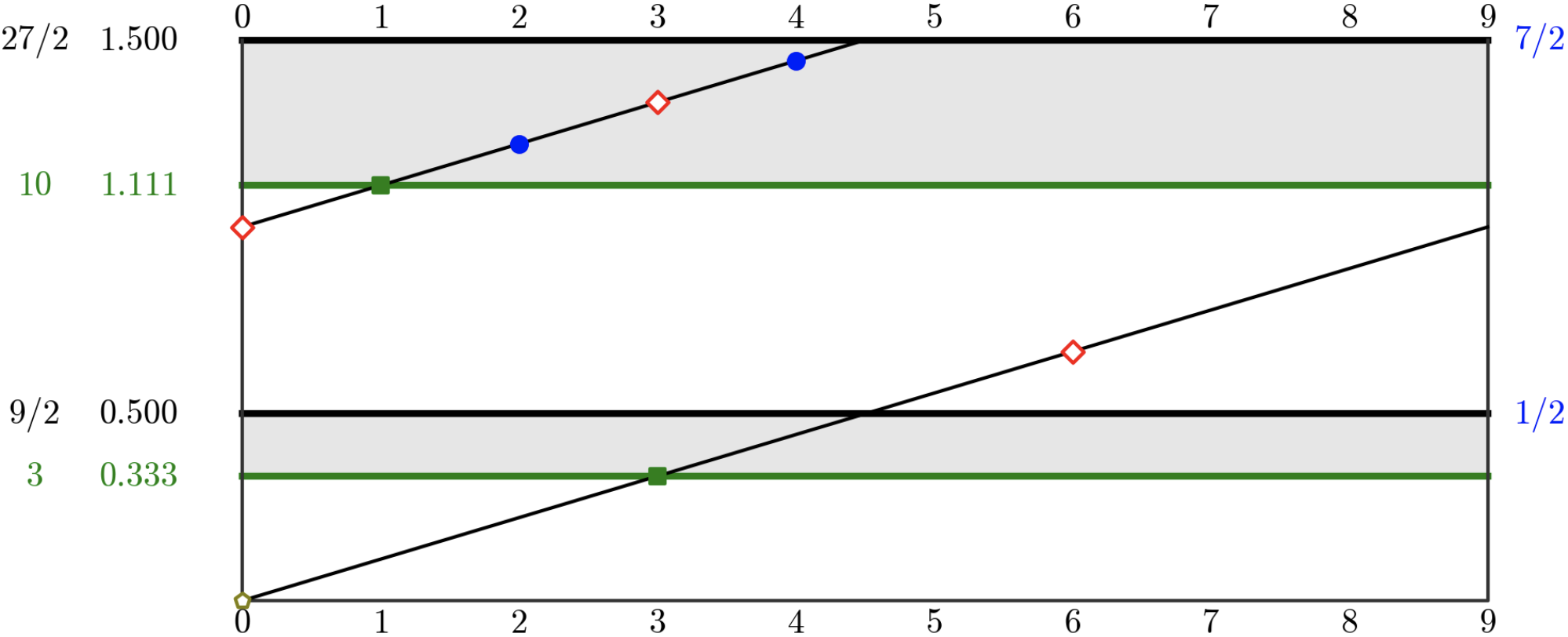}
\end{center}
\caption{\label{fig examp9} The Eisenstein diagram of the introductory nonic $3$-adic Herbrand invariant with
 \gr{means $\langle \frac{1}{3},\frac{10}{9}\rangle$}, 
 slopes $[\frac{1}{2},\frac{3}{2}]$, \bl{rams $(\frac{1}{2},\frac{7}{2})$}, and generic polynomial \eqref{eq nonicgenpoly}.
Various general facts are relatively easy to see here because the two bands
do not overlap.} 
\end{figure}

Figures~\ref{fig examp9} and \ref{fig examp8} each show the Eisenstein diagram 
for a family discussed previously, as described in their captions.  
 Most of our discussion of Eisenstein diagrams is
supported by one or both of these figures. We will give links to the database for illustrations of other phenomena.

\subsubsection*{A spiral with points representing terms} We think of the rectangle $R$ as 
a cut and unrolled cylinder, with each horizontal line segment of constant level $s$ coming 
from a circle.  We draw the spiral that starts at 
$(0,0)$ and goes up with slope $1/e$.   This spiral
wraps every time it meets the right edge of $R$.
The part of the spiral that goes from 
$(0,j)$ to $(e,j+1)$ is called the $j^{\rm th}$ ramp.   

We place conditions on the term $\pi f_\sigma \pi^j x^i$ by referring
to the unique point $P_\sigma$ on ramp $j$ with horizontal
coordinate $i$.   Equivalently $P_\sigma$ is the unique point
on the spiral with $s = \sigma/e$.  Thus $P_0 = (0,0)$ is the 
starting point of the spiral, and as one goes up one
encounters the points $P_1$, \dots, $P_{\lfloor e s_w \rfloor}$ in
order.  

\subsubsection*{Bands} 

For $k=1$, \dots, $w$,  the band $B_k$ is defined to be the set of 
points in $R$ satisfying $m_k  \leq s < s_k$.  
 As $k$ increases, the $m_k$ strictly increase
and so the bottom edges of bands go up.  However the
$s_k$ only weakly increase.   So the top edges of the
bands $B_k$ for $k$ in a segment of common $s_k$ agree, 
even as the widths of these bands successively decrease by a factor of
$p$, as in \family{2.1.16.30a}, the sample family summarized in Table~\ref{tab sixteentab}. 
It will also be 
convenient to let $B_0$ denote the bottom edge
of the rectangle, i.e.\ the points with $s=0$.   

We shade bands by gray and indicate
overlaps by darkening the gray.  We color bottom edges
of bands green and the bounding top edges black.  
In the rare cases where an upper edge agrees with a lower
edge we dash this boundary using black and green, as in 
\family{2.1.16.79a}, the family highlighted in \S\ref{subsec globalization}. 
These conventions assist in visually identifying the bands,
even when they overlap.

\subsubsection*{Index and types of points} We say that the {\em index}
 of an integer $i \in [0,e)$ is $k = w-\mbox{ord}_p(\gcd(i,p^w)) \in \{0,1,\dots,w\}$.   We partition
 the points $P_\sigma$ into five types, using colors and shapes 
 to distinguish the types.  
 
 The point $P_0 = (0,0)$ plays a special role.  We call it the {\em $D$-point} and we mark it by an olive pentagon. 
A $D$-point is {\em critical} and 
drawn solid if $\gcd(\epsilon,q-1)>1$.
Otherwise it is {\em negligible} and drawn hollow.  
 Here, like with the red diamonds below,
 the solidness indicates complications in the process 
 of choosing unique distinguished polynomials for fields $L/K$ in the family $I/K$.
 
  For $\sigma>0$, we use 
 the band $B_k$ to classify the points $P_\sigma$ of index $k$ as follows.  
\begin{itemize}
\item {\em $Z$-points} are points beneath $B_k$.  Their 
associated color is {\em clear,} meaning we don't draw them.   
\item  {\em $A$-points} are the unique points $P_\sigma$
which are at the bottom edge of geometric bands $B_k$ for which
$k$ is final in its segment, meaning that either $k=w$ or $s_{k+1}>s_k$.   
We draw them as 
{\em solid green 
squares}. 

\item {\em $B$-points} are points in their governing band $B_k$
which are not $A$-points.   They are represented by {\em solid blue disks.} 

\item {\em $C$-points} are points above their governing band.  
They are represented by {\em red diamonds.}   A $C$-point is {\em critical} and 
drawn solid if it is on the top border of any band.  Otherwise
it is {\em negligible} and drawn hollow.  
\end{itemize}
The infinite set $\Eisen(I/K)$ now has the following explicit description.  It is the subset
of $\Eisen(e/K)$ where $f_\sigma=0$ for $Z$-points and $f_\sigma \neq 0$ for
$A$-points.  The key finite set $\EisenII(I/K)$ is the subset of $\Eisen(I/K)$ where 
$f_0=1$ if the $D$-point is negligible and  
$f_\sigma=0$ for negligible $C$-points.    The fact that $\EisenII(I/K)$ still represents all 
extensions in $I/K$ is not at all obvious.  The 
proof of a considerably stronger statement is by Monge's reduction algorithm \cite[\S 2]{Mo}.  Basic
aspects of this theory will be described in \S\ref{subsec reduction} below, 
and the full theory is described in a computational
context in \cite{polred}.

\cmmt{
A point is critical if and only if the corresponding relative
extension is rigid in the sense of having no non-trivial automorphisms after an unramified extension of $K$.  
For the $D$-point the relative
extension is $L_0/L_{\rm un}$ and so the associated
inertia group is cyclic of order $\epsilon$.  
The automorphism group has order $\epsilon$ which is
why $\epsilon>1$ is the 
For the $A$-points appearing on the top
edge of $\rho$ bands, the extension $L_k/L_{k-\rho}$ has
associated inertia group of the form $C_p^\rho \rtimes C_d$ 
with $d$ the denominator of the ram $r_k$.  So the extension
has such automorphisms if and only if $d=1$, i.e.\ if and
only if $t_k$ is integral.  
}

\subsection{Generic polynomials and numerics} 
\label{subsec generic} 
\begin{figure}[t]
\begin{center}
\includegraphics[width=4.5in]{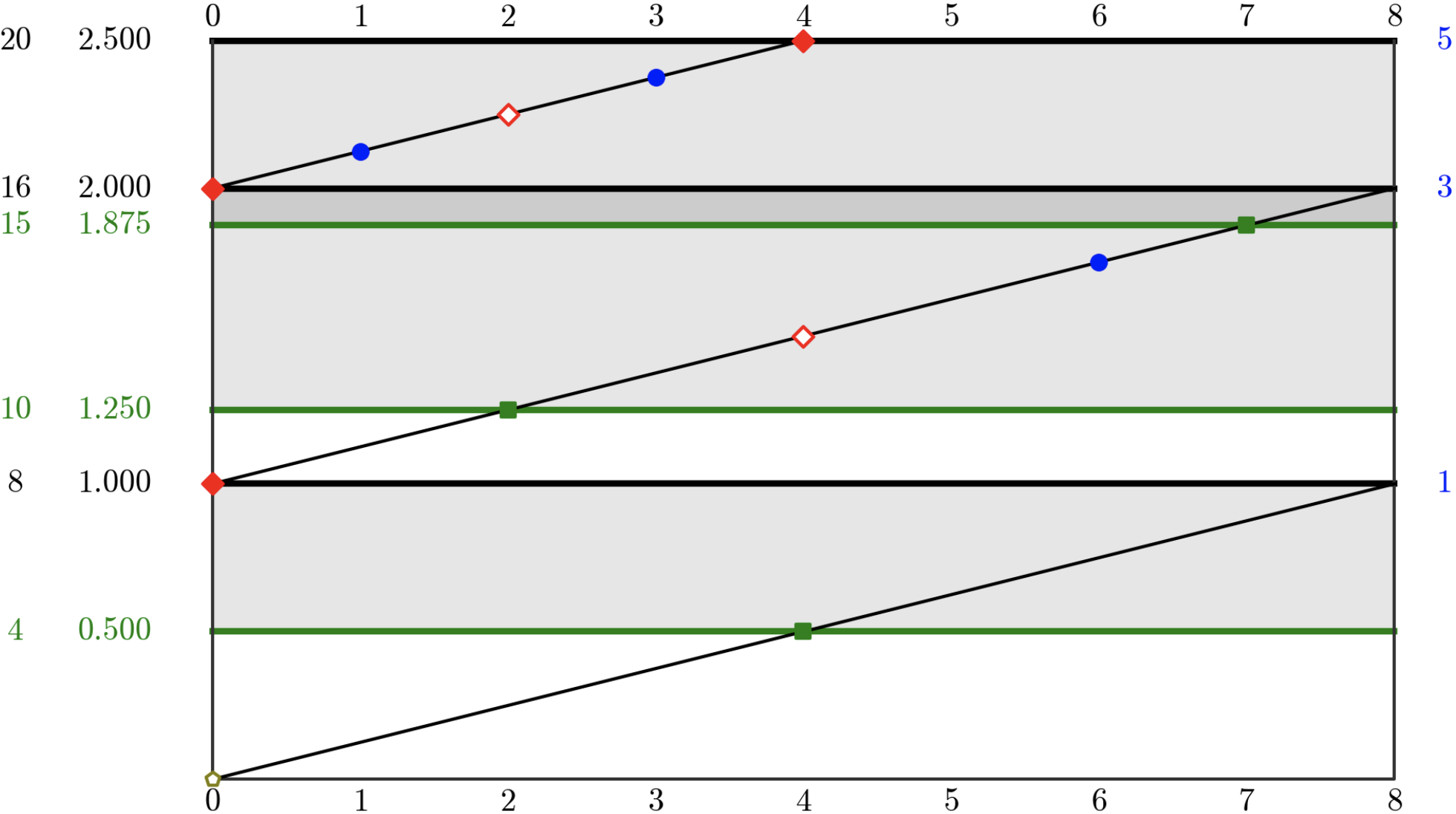}
\end{center}
\caption{\label{fig examp8} Eisenstein diagram for the octic $2$-adic 
Herbrand invariant with \gr{means
 $m = \langle  1,\frac{5}{2},\frac{15}{8} \rangle$},
 slopes ${s=[1,2,\frac{5}{2}]}$, and \bl{rams $r = (1,3,5)$}.   
 }
\end{figure}

The generic polynomial for a given 
Herbrand invariant has the form $\left( \pi \sum_\sigma f_\sigma \pi^j x^i \right) + x^e$.  
Here $\sigma$ runs over all nonnegative integers for which $P_\sigma$ is drawn,
but not negligible.  To make structure more evident, we replace $f_\sigma$ by 
$d_\sigma$, $a_\sigma$, $b_\sigma$, or $c_\sigma$ according to whether
$P_\sigma$ is a $D$-point, $A$-point, $B$-point, or $C$-point.    So when these coefficients
run independently over $\tilde{\kappa}$, except for the inequalities 
$d_0, a_\sigma \neq 0$, one
gets the set $\EisenII(I/K)$.    As an example beyond \eqref{eq nonicgenpoly}, the generic 
polynomial corresponding to the $2$-adic Herbrand invariant $[1,2,2.5] = (1,3,5) = \langle 0.5,1.25,1.875 \rangle$ 
of Figures~\ref{fig slopepolygon}, \ref{fig herbrandphi} and \ref{fig examp8} is 
 \begin{eqnarray}
\nonumber \lefteqn{f(a_4,a_{10},a_{15},b_{14},b_{17},b_{19},c_8,c_{16},c_{20}; \pi;  x)  =  (\pi + \pi^2 c_8 + \pi^3 c_{16}) +  \pi^3 b_{17} x + } \\
\label{eq octicgenpoly}&&  
   \pi^2 a_{10} x^2 + \pi^3 b_{19} x^3+   (\pi a_4 + \pi^3 c_{20}) x^4 +
  \pi^2 b_{14} x^6+
 \pi^2 a_{15} x^7+ 
 x^8.
 \end{eqnarray}
 Other examples are given in the next sections.
 
Let $\delta \in \{0,1\}$ be the number of critical $D$-points, and let $\alpha$, $\beta$, 
and $\gamma$ be the number of $A$-points, $B$-points, and critical $C$-points respectively.
 Clearly $|\EisenII(I/K)| = (q-1)^{\delta+\alpha} q^{\beta+\gamma}$.  Monge
 reduction (see \S\ref{subsec reduction}) says that $|\EisenII(L/K)| = (q-1)^{\delta} q^{\gamma}/|\Aut(L/K)|$.
 Accordingly, one gets a mass formula 
  \begin{equation}
 \label{eq massformula}
 M(I) := \sum_{L/K \in I/K} \frac{1}{|\Aut(L/K)|}= (q-1)^\alpha q^\beta.
 \end{equation}
 In the rigid case $\delta = \gamma=0$ where $|\Aut(L/K)|=1$ is forced,
 \eqref{eq massformula} becomes a cardinality formula. 
   As a non-rigid example, from \family{2.1.8.22d} there are 
 $32$ fields in the family \eqref{eq octicgenpoly}, with \eqref{eq massformula} becoming $8 (\frac{1}{8}) + 20 (\frac{1}{4}) + 4 (\frac{1}{2}) = 8$.  

 The numbers $M(I)$ in Figure~\ref{fig families} come from \eqref{eq massformula} with $q$ set equal to either $2$ or $3$.
 A clarifying check on various numbers is 
 to view the space $\Eisen(e/K)$ of degree $e$ Eisenstein polynomials over $K$ as an infinite product $\F_q^\times \times \F_q \times \F_q \times \cdots$, with the successive
 factors corresponding to the variables
 $f_0$, $f_1$, $f_2$, \dots.  Giving each
 factor its uniform probability measure 
 turns $\Eisen(e/K)$ into a probability 
 space.  Serre proved in \cite[Thm.\,2]{serrepaper}
 that the chance that a random polynomial
 in $\Eisen(e/K)$ is in $\Eisen(I/K)$ 
 is the Serre mass $SM(I) = M(I)/q^{e m_w}$.
 Figure~\ref{fig families} lets one see the
 terms in the resulting formula 
 $\sum_{I/K} SM(I) = 1$ for three
 different $(p,e_K)$.  Taking the easiest
 case $(2,1)$ as an example, there are
 seven families, and one can use $em_2 = 2r_1+r_2$. The sum is the dot 
 product $(1,2,2,2,4,4,8) \cdot 2^{-(1,3,5,5,6,7,8)}$ and it is indeed $1$.

\subsection{Slopes, means, and rams on Eisenstein diagrams}
\label{subsec slopesmeansrams}  When drawing
 Eisenstein diagrams,
both here and on the family pages of the database, we indicate the slopes, means, and rams
in appropriate places.   There is no need to 
tick the vertical axis at the left of the diagram in the traditional way, because the left endpoint
of the $j^{\rm th}$ ramp is the integer $j$.  Instead we give decimal
approximations to lower edges $m_k$ and upper edges $s_k$.  
In the leftmost column, we also give the scaled versions $em_k$ and $es_k$.
These scaled versions make the subscripts on the marked points $P_\sigma$ more immediately 
identifiable.  

The ram $r_k$ is printed to the immediate right of the upper boundary
of the band $B_k$.   It gives two point counts as follows.  
First, keeping in mind that $B_k$ includes its lower boundary but not its 
upper boundary, the number of drawn points in $B_k$ of index
at most $k$ is $\lceil r_k \rceil$.  
Second, the total number of red diamonds on or below the top edge of $B_k$ 
is $\lfloor r_k \rfloor$.   

We have just given two interpretations of rams, but there are two easier ones.  For these further 
interpretations, 
it is useful to scale rams to
\begin{align}
\label{eq scaled} 
\mbox{\em small rams } r_k^* & = \frac{(p-1)r_k}{\epsilon p^k} &
\mbox{and {\em tiny rams} } \;\; r_k' & = \frac{r_k}{\epsilon p^k}.
\end{align}
Then
\begin{align}
\label{eq twotransforms}   r_k^* & = m_k-m_{k-1} & \mbox{and } \;\;  r'_k & = s_k-m_k.
\end{align}
Together with the standing convention $s_0=m_0=r_0=0$, these simple equations
give the complete relations between slopes, means, and rams.  So the tiny ram $r_k'$
is immediately visible as the vertical width of the $k^{\rm th}$ band.

The set of Herbrand invariants $I = (r_1,\dots,r_w)_\epsilon^f$ compatible with a 
given ground field $K$ has an appealing geometric description if one uses small rams.  
Namely $(r_1^*,\dots,r_w^*)$ must be in 
the cube $(0,e_K]^w$, with all $r_k$ geometric
if it is in the interior $(0,e_K)^w$.     
The hook $H_j$ of \S\ref{subsec classification} indicates a part of the boundary of $[0,j]^2$.  In our Eisenstein 
diagrams, arithmetic bands, meaning the ones
where $r_k^* = (p-1) r_k'$ is an integer, as discussed in \S\ref{subsec classification}, are marked by a black segment at their right end. Thus the bottom two bands in  \family{3.1.27.99c} are marked, while the top one is not.   

\cmmt{
The geometric vs.\ arithmetic distinction of \S\ref{subsec classification} for rams is clearer when
transported to small rams.  Namely if $I/K$ is a family, then one 
has always $t_k^* \leq e_K$.  The ram is geometric if there is strict inequality 
and arithmetic if there is equality.   Similarly, whether or not a given
Herbrand invariant $I$ occurs over a given ground field $K$ becomes clearer as
follows.   Define the {\em size} of a wild Herbrand invariant $I$ to be $t^* = \max t_k^*$.
If $t^*$ is integral then $I$ occurs exactly if $t^*=e_K$.  If 
$t^*$ is nonintegral than $I$ occurs exactly if $t^*<e_K$. }

\subsection{The Herbrand function and Monge reduction}
\label{subsec reduction}
Of the four versions \eqref{eq herbrand}-\eqref{eq herbrandrams} of totally ramified Herbrand invariants, the only 
one that is not immediately evident from an Eisenstein diagram is the first one, involving 
the classical Herbrand function $\phi$.  To read $\phi$ off from an Eisenstein diagram,
one can think of a particle starting at time $r=0$ at the point $P_0=(0,0)$ and moving up the spiral. 
The particle starts at a speed of $p^w$ steps per second, where moving from $P_\sigma$ to $P_{\sigma+1}$ counts
as a step.  
Every time the particle crosses the top edge of a band it decreases its speed by a factor 
of $p$. This is a natural definition of speed, because the particle encounters $C$-points exactly at
positive integral times. At any time $r \in \R_{\geq 0}$,
 the particle has traveled some number of ramps $s \in \R_{\geq 0}$. 
Then, as in Figure~\ref{fig herbrandphi}, $s=\phi(r)$.    

  The Monge reduction algorithm iteratively simplifies a given 
  Eisenstein polynomial without changing the field it defines. 
  Its general nature can be understood
  in terms of the moving particle.  At the initial $D$-point
  a change of variables is made trying to make $f_0=1$.  
  At each subsequent $C$-point $P_\sigma$ a change of 
  variables is made trying to make $f_\sigma = 0$.  These
  changes to $f_\sigma$ also change some of the $f_{\sigma'}$ for 
  $\sigma'>\sigma$.   When the particle leaves the 
  drawn window, the complicated process can be stopped. 
  Instead one can just turn all the $f_\sigma$ with $\sigma>e s_w$ 
  to $0$, as this change does not affect the field defined, 
  by an effective version of Krasner's Lemma.  
  
  The reduction process is completely successful at negligible points
  but only partially successful at critical points.  If one simply does not 
  make the coordinate change at critical points, then one gets a
  surjection from the infinite set $\Eisen(L/K)$ to the $(q-1)^\delta q^\gamma$-element
  set $\EisenII(L/K)$.   If one makes the coordinate 
  changes at the critical points as well, then the ambiguity 
  from a critical $D$-point reduces
  from $q-1$ to $\gcd(q-1,\epsilon)$, and the ambiguity from a
  critical $C$-point at the top of $\rho$ bands reduces from $q$ to 
  at most $\gcd(q,p^\rho)$. Multiplying these bounds together gives an 
  ambiguity constant $\Amb(I/K)$.  
  It depends on $K$ only through $q$ and is a divisor of the 
  integer $\Amb(I)$ introduced in \S\ref{subsec rigidity}. 
  The full reduction algorithm gives a surjection from $\Eisen(L/K)$ to 
  a subset of $\EisenII(L/K)$ of
  size a divisor of $\Amb(I/K)/|\Aut(L/K)|$.  The sequel \cite{polred} implements the 
  full algorithm, and moreover deterministically chooses a distinguished
  polynomial from the set of outputs.

\section{Sample families}
\label{sec sample}
    This section is aimed at facilitating the reader's exploration of the database.  
It summarizes the internal structure of several families, emphasizing topics
which support the more theoretical discussions of the next section.

\subsection{Partitions of the introductory family $[\frac{1}{2},\frac{3}{2}]/\Q_3$}
\label{subsec packets}
All absolute families are naturally partitioned in two related ways, into subfamilies
and into packets.  The subfamilies are determined by using the coefficients 
corresponding to points at the bottoms of bands; these are $d_0$, the $a_\sigma$,
and sometimes also some $b_\sigma$, as discussed in
\cite{polred} in the context of residual polynomials.  Our
naming convention for fields incorporates subfamilies,
because subfamilies are both elementary and important.

\begin{table}[htb]
 \[
{\renewcommand{\arraystretch}{1.2}
\renewcommand{\arraycolsep}{3.5pt}
\begin{array}{cccc|cc}
 a_3 & a_{10} & b_{11} &  b_{13}   &  \mbox{HSC} & G   \\
 \hline
1& 1 & 1 &  0 &  [\;]_2 & 9T5 \\
1& 1 & 0,2 &  0,1,2 &    [1]_2 & 9T10 \\
1& 1 & 1&   1,2 &   [\;]^3_2 & 9T11 \\
\hline
1& 2 & 1 &  0,1,2 &  [1]^2_2 & 9T18 \\
1& 2 & 2 & 0,1,2 &   [\frac{1}{2},\frac{2}{3}]_2^2 & 9T24 \\ 
1& 2 & 0 &  0,1,2 &  [\frac{1}{2},1]^2_2 & 9T24 \\
\end{array}
}
\;\;\;\;\;\;\;
{\renewcommand{\arraystretch}{1.2}
\renewcommand{\arraycolsep}{3.5pt}
\begin{array}{cccc|cc}
 a_3 & a_{10} & b_{11} &  b_{13}   &  \mbox{HSC} & G   \\
 \hline
2 & 1 & 2 &  2 &   [\;]_2 & 9T3 \\
2 & 1 & 2 &  0,1 &   [\;]^3_2 & 9T10 \\
2&  1 & 0 &  0,1,2 &   [1]_2 & 9T10 \\
2 & 1 & 1 &  0,1,2 &  [1]_2 & 9T11 \\
\hline
2 & 2 & 1 &  0,1,2 &   [\;]^2_2 & 9T8 \\
2 &  2 & 0,2 &  0,1,2 &     [\frac{1}{2},1]^2_2 & 9T24 \\
\end{array}
}
\]
\caption{\label{tab nonicpackets} Partitions of the
introductory family $[\frac{1}{2},\frac{3}{2}]/\Q_3 = \family{3.1.9.18b}$ into
two subfamilies and ten packets.  
}
\end{table}
Giving a name to a concept introduced in \S\ref{subsec new},  we say that two fields
are in the same {\em packet} if they have the same Galois 
group $G$ and the same hidden slope content $HSC$.
The database has $(G,HSC)$ for all its fields in degree $\leq 15$.
However since $G$ and particularly $HSC$ can be hard to compute in higher 
degrees, packets are not incorporated into our labeling scheme.

As the sample family of this subsection, we reconsider the
family $[\frac{1}{2},\frac{3}{2}]/\Q_3$ of the introduction. 
Its generic polynomial from \eqref{eq nonicgenpoly} is 
$3 + 9 a_{10} x  + 9 b_{11} x^2 +  3 a_3  x^3  +  9 b_{13} x^4  + x^9$.
Table~\ref{tab nonicpackets} breaks the $36$ fields into two subfamilies
of eighteen fields each, according to the value of $a_{10}$.
The canonical cubic subfield is defined by the Eisenstein polynomial $3 + 3 a_3 y + y^3$,
so the relative families discussed at the end of \S\ref{subsec relative} correspond to 
the left and right halves of Table~\ref{tab nonicpackets}.  
In a more complicated way, the table breaks the family
into its ten packets.  Thus the packet discussed in \S\ref{subsec new}
is given in the second from bottom line in the left half.

\subsection{Easy linear packets in $[\frac{3}{2}]_5/\Q_3$} 
\label{subsec easylinear}
\cmmt{To give a reader a better feel for the level of complication
one can expect in the packet structure of families,
consider the collection of families over $\Q_p$ 
having $f=1$, any possible $\epsilon$, and exactly one wild slope 
$s=k/(\epsilon(p-1))$.     So the
 indexing triple $(p,\epsilon,k)$ consists of an arbitrary
 prime number $p$, and arbitrary positive integer $\epsilon$
 prime to $p$, and, by \S\ref{subsec classification}, an integer $k$ in 
 $[1,\dots,p\epsilon]$ which is either prime to $p$ or 
 $p \epsilon$ itself.}

The family $[\frac{3}{2}]_5/\Q_3$ with LMFDB label \family{3.1.15.29a} has a very 
simple packet structure as follows. 
The generic polynomial is 
$3+9(b_{16} x+ b_{17} x^2 + b_{19} x^4 + b_{20} x^5 + b_{22} x^7) + x^{15}$.  
Specializing via $b_\sigma \in \{0,1,2\}$ gives $3^5=243$ polynomials
bijectively representing the $243$ fields of the family.   If the list 
$(b_{16}$, $b_{17}$, $b_{19}$, $b_{22}$) starts
with exactly $0$, $1$, $2$, or $3$ zeros, then the hidden slope
content is $[j/10,j/10,j/10,j/10]_2^4$ for $j = 13$, $11$, $7$, $1$, and the
Galois group is $15T64=C_3^4: (S_3 \times F_5)$.   If the list
is simply $(0,0,0,0)$, then the hidden slope content is $[\;]_2^4$, and the 
Galois group is $15T11 = S_3 \times F_5$.  The fact that the coordinates
$b_\sigma$ render the packet structure transparent is an example of the
``important bonus'' mentioned 
after \eqref{eq nonicgenpoly}.

\subsection{Complicated linear packets in $[2]_7/\Q_2$} 
\label{subsec hardlinear}
 The
family $[2]_7/\Q_2$ with LMFDB label \family{2.1.14.27a} has
generic polynomial 
\[
(2+4b_{14}+8c_{28}) + 4\left( b_{15} x +  b_{17} x^3 +  b_{19} x^5 + b_{21} x^7 +  b_{23} x^9 +  b_{25} x^{11} +  b_{27} x^{13}\right) + x^{14}.
\]
There are sixteen packets $(G,HSC)$, with the hidden slope content $HSC$
always determining the Galois group $G$.  The possibilities
 for the wild slopes are indexed by the set $\{13,5,3,\emptyset\} \times \{11,9,1,\emptyset\}$.
 Here an index $j$ generically contributes 
  $[j/7,j/7,j/7]$ but $\emptyset$ contributes the
 empty list.  Inspecting the database shows that the list of wild slopes does not depend on
 $b_{21}$ and $c_{28}$.    Table~\ref{tab fourteenpackets} shows the dependence on the remaining coefficients.
 In the displayed vectors, each $\star$, $v$, or $d$ can 
be independently $0$ or $1$, except that the $d$'s must sum to an odd number
and the $v$'s to an even number.  

 \begin{table}[htb]
 \[
 {
 \begin{array}{c|cccc}
 & 11 & 9 & 1 & \emptyset  \\
 \hline
 13 & (1,0,\star,\star,\star,\star) & (1,1,1,\star,\star,\star) & (1,1,0,d,d,d) & (1,1,0,v,v,v)  \\
 5 &   (0,1,\star,\star,0,\star) & (0,0,1,1,\star,\star) & (0,0,0,1,\star,1) & (0,0,0,1,\star,0)   \\
 3 &   (0,1,d,1,d,\star) & (0,0,1,0,0,\star) & (0,0,0,0,1,1) & (0,0,0,0,1,0)   \\
 \emptyset & (0,1,v,1,v,\star) & (0,0,1,0,0,\star) & (0,0,0, 0,0,1) & (0,0,0,0,0,0)  \\
 \end{array} 
 }
 \]
 \caption{\label{tab fourteenpackets} 
 The hidden slope stratification of the family $[2]_7/\Q_2$, with vectors 
$(b_{15},b_{17},b_{19},b_{23},b_{25},b_{27})$ indicating coefficients.}
 \end{table}

\subsection{Families $[1,\dots,1]/\Q_p$}  
\label{subsec unisloped} 
Families where $\epsilon=1$ and there is just a single wild slope $s$ repeated
$w$ times are easier, because there are no hidden wild slopes and 
Galois groups can be completely determined as in Theorem~8.2 of
\cite{greve-pauli}.
The case where $s=1$ is particularly interesting for several reasons.
This subsection focuses on the families $[1,\dots,1]/\Q_p$, there being
one for each prime power $p^w$, using the particular family
$[1,1,1,1]/\Q_2$ as an example.

The generic polynomial for the $2$-adic Herbrand invariant $[1,1,1,1]$ is 
\begin{equation}
\label{eq sixteengenpoly}
\pi \left((1+\pi c_{16}) + b_8 x^8 + b_{12} x^{12} + b_{14} x^{14} + a_{15} x^{15}\right) + x^{16}.
\end{equation}
A field $L/\Q_2$ in $[1,1,1,1]/\Q_2$ has one or two representing polynomials according to whether
$b_8+b_{12}+b_{14}$ is even or odd.  In the latter case, the two polynomials differ only in
the coefficient $c_{16}$.  The polynomial with $c_{16}=0$ is our choice of distinguished polynomial.  
\begin{table}[htb]
{
\[
{\renewcommand{\arraystretch}{.3} 
\renewcommand{\arraycolsep}{2pt}
 \begin{array}{ccccccccccccc}
  \;\; \ell.i \; & \,b_8 \,& b_{12}  &  b_{14} &  c_{16} & 
                     \mbox{Associated polynomial} &\;\; \; \; a \; \; \;\;& \; \; \;\; u \; \; \; \; & \; [j_0, & j_1, & j_2, & j_3, & j_4] \\
 \hline
                 \ffield{2.1.16.30a}{1.1} & 0 & 0 & 0 & 0 & (y+1)^4 & 2 & 4 & [1, & 3, & 7, & 15, &
                     30] \\
                   \ffield{2.1.16.30a}{1.2} & 0 & 0 & 0 & 1 & (y+1)^4 & 2 & 4 &  [1,& 3, & 7, & 15, &
                     32] \\
                     &\\
                   \ffield{2.1.16.30a}{2.1} & 0 & 0 & 1 & 0 & y^4+y+1 & 1 & \; 15 & [1, & 3, & 7, & 14,
                     & 31] \\
                     &\\
                   \ffield{2.1.16.30a}{3.1} & 0 & 1 & 0 & 0 & \left(y^2+y+1\right)^2 & 1 & 6  & [1,
                     & 3, & 6, & 12, & 31] \\
                     &\\
                   \ffield{2.1.16.30a}{4.1} & 0 & 1 & 1 & 0 & (y+1) \left(y^3+y^2+1\right) & 2 & 7 
                     & [1, & 3, & 6 ,& 15, & 30] \\
                   \ffield{2.1.16.30a}{4.2} & 0 & 1 & 1 & 1 & (y+1) \left(y^3+y^2+1\right) & 2 & 7 
                     & [1, & 3, & 6, & 15, & 32] \\
                     &\\
                   \ffield{2.1.16.30a}{5.1} & 1 & 0 & 0 & 0 & y^4+y^3+1 & 1 & 15 &  [1, & 2, & 4, & 8,&
                     31] \\
                     &\\
                  \ffield{2.1.16.30a}{6.1} & 1 & 0 & 1 & 0 & (y+1)^2 \left(y^2+y+1\right) & 2 & 6 
                     & [1, & 2, & 4, & 15, & 30] \\
                   \ffield{2.1.16.30a}{6.2} & 1 & 0 & 1 & 1 & (y+1)^2 \left(y^2+y+1\right) & 2 & 6 
                     & [1, & 2, & 4, & 15, & 32] \\
                     &\\
                   \ffield{2.1.16.30a}{7.1} & 1 & 1 & 0 & 0 & (y+1) \left(y^3+y+1\right) & 2 & 7 
                     & [1, & 2, & 7, & 15, & 30] \\
                   \ffield{2.1.16.30a}{7.2} & 1 & 1 & 0 & 1 & (y+1) \left(y^3+y+1\right) & 2 & 7 
                     & [1, & 2, & 7, & 15, & 32] \\
                     &\\
                   \ffield{2.1.16.30a}{8.1} & 1 & 1 & 1 & 0 & \;\; \;\;y^4+y^3+y^2+y+1 \;\; \;\;& 1 & 5 & [1, & 2, & 7,
                     & 14, & 31] \\
                  \end{array}
                  }
   \]
   }
   \caption{\label{tab sixteentab} Information on the $12$ fields \texttt{2.1.16.30a}$\ell.i$ in the family \family{2.1.16.30a}}
  \end{table}
  Table~\ref{tab sixteentab} presents information directly available on the family page for $[1,1,1,1]/\Q_2$, in a form slightly 
  modified to support the discussion here.    In the first column, $\ell$ indexes the subfamily and
  $i$ indexes the field within the subfamily.   The jump sets $[j_0,j_1,j_2,j_3,j_4]$ are discussed
  in general in \S\ref{subsec jumpsets}.  Commonly, a family gives rise to just a very few jump sets, often just
  one.  This family, and conjecturally all the $[1,\dots,1]/\Q_2$, have the unusual feature
  that the jump set determines the field.   
  
  The general case $[1,\dots,1]/\Q_p$ has a generic polynomial of a form similar to \eqref{eq sixteengenpoly}.  
 Write $v_w=a_{p^w-1}$ and $v_j = b_{p^w-p^{w-j}}$ for $j=1$, \dots, $w-1$.   Then the parameters in the generic polynomial are $v_1$, \dots, $v_w$ and also $c_{p^w}$.   Let
  $L/\Q_p$ be the field defined by the parameters $(v_1,\dots,v_w;c_{p^w})$.  
 Let $g$ be an  element of $GL_w(\F_p)$ with characteristic
 polynomial 
\begin{equation}
\label{residual} f(y) = y^w + \sum_{j=0}^{w-1}  v_{w-j} y^j.  
\end{equation}
Then the Galois group of $L^{\rm gal}/\Q_p$ is the semidirect
product $\F_p^w \rtimes \langle g \rangle$  \cite[Theorem 8.2]{greve-pauli}.  Here the wild inertial group 
is $\F_p^w$, the tame quotient of inertia has order one, and the unramified
quotient is $\langle g \rangle$.  Thus the hidden slope content is simply $[ \; ]^u$,
where $u$ is the order of $g$.    
For the case $p^w = 2^4$, Table~\ref{tab sixteentab} gives these
polynomials in factored form.  The residual polynomials given in
the database can be obtained from $f(y)$ by replacing each $y^i$ with
$z^{2^i}$ and dividing by $z$.
The database gives the Galois groups in the usual way,
from the smallest group $16T166= \F_2^4 \rtimes C_4 = C_2 \wr C_4$ to the largest group 
$16T447 = \F_2^4 \rtimes C_{15} = F_{16}$.

The various phenomena discussed in the example of $2^4$ generalize to $p^w$.  
A field $L/\Q_p$ has either one or $p$ representing polynomials according to 
whether $f(1)$ is different from or equal to zero in $\F_p$.  In the latter case, 
the polynomials again differ only in the coefficient $c_{p^w}$ and again
$c_{p^w}=0$ gives our choice of distinguished polynomial.  So, as illustrated
on the table for $2^4$, the automorphism number $a = |\Aut(L/\Q_p)|$ is
$p$ if $f(y)$ has a factor $(y-1)$ and $1$ otherwise.  Going further, the subfields of 
$L/\Q_p$ are in natural bijection with the
factors of $f(y)$ in $\F_p[y]$, with the field
coming from a degree $k$ factor having
degree $p^k$.   So there are all together $\prod (m_j + 1)$ subfields,  where $\prod_j f_j(y)^{m_j}$ is the factorization
of $f(y)$ into irreducibles.  

\cmmt{
Now consider changing the ground field from $\Q_p$ to $\F_p((t))$. To agree with
standard conventions in the positive characteristic literature, make the replacements 
$x \rightarrow 1/x$ and $\pi \mapsto -1/t$ in the generic polynomial, and 
partially clear denominators by multiplying 
through by $t x^{p^w}$. Setting the polynomial equal
to $0$ and separating the terms involving $t$ from 
those involving $x$ gives the equation
\kk{I think the last term should be multiplied by $x^{p^w}$.}
\begin{equation}
\label{eq cover}
x^{p^w} + \sum_{j=0}^{w-1} u_j x^{p^j} =  t+ \frac{c_{p^w}}{t}.
\end{equation}
We have set up this equation as defining an extension of the
local field $\F_p((t))$.  But it also defines an extension
of the global field $\F_p(t)$.  

When $c_{p^w}=0$, as it is often in our context, Equation~\ref{eq cover}
is particularly appealing.  Then the equation 
can be interpreted as a map from the projective $x$-line to the
projective $t$-line.   The derivative of the left side is the constant $u_0 = a_{p^w-1} \neq 0$,
so this map is ramified only over $\infty$.  The left side is an additive polynomial,
so the Galois group is a subgroup of the affine group 
$\mbox{Aff}_w(\F_p) = \F_p^w \rtimes GL_w(\F_p)$.   Because the 
coefficients are constant, the Galois group in fact agrees
with the decomposition group at $\infty$, described as 
$\F_p^w \rtimes \langle g \rangle$ above.   We have
just illustrated the concept of {\em globalization}, 
which we will discuss more in \S\ref{subsec globalization}}

\section{Theoretical connections}
\label{sec theoretical}
Consider finite extensions $K$ of 
$\Q_p$ inside of a fixed separable closure $\Qpsep$.  A natural goal, which seems a long way off or perhaps not obtainable at all, is to find a group-theoretical description of each
absolute Galois group $\Gal(\Qpsep/K)$ together
with its descending filtration by the higher ramification groups
$\Gal(\Qpsep/K)^s$.  There are however many deep theorems
towards this goal.  In this section, we describe ways in which the 
database interacts with these theorems, rendering them more
explicit.

\subsection{Extensions with a given Galois group} 
\label{subsec givenGalois} An overview of many powerful and explicit 
results on absolute Galois groups is given in 
\cite[VII\S5]{neukirch-schmidt-wingberg}.
Highlights are as follows. Let 
$K^{\rm nil}$ be the maximal extension of $K$ for which the Galois group is a pro-$p$-group; here ${\rm nil}$ stands for nilpotent.  
Shafarevich proved in the 1940s that 
$\Gal(K^{\rm nil}/K)$ is free on $[K:\Q_p]+1$ generators
if $K$ does not contain $p^{\rm th}$ roots of unity.  Demushkin proved in the 1950s that it can always be presented with $[K:\Q_p]+2$ generators and one explicit relation, with an example
being 
\begin{equation}
\label{eq G2}
\Gal(\Q_2^{\rm nil}/\Q_2) = \langle x,y,z |x^2 y^4 x^{-1} y^{-1} x y \rangle.
\end{equation}
In the 1980s, Jannsen and Wingberg 
gave a complete description of  
$\Gal(\Qpsep/K)$ for 
$p$ odd, and then Diekert did 
the same for $p=2$, assuming
$K$ contains fourth roots of unity.  
All these results are silent on the
ramification filtration.

\cmmt{
From the 1980s,
one has a complete description of $\Gal(\Qpsep/K)$ 
by Jannsen and Wingberg for $p$ odd and by Diekert for $p=2$ if $K$ contains 
fourth roots of unity.  As a similar result not covered by this general theorem, \jg{we could specify which one} let 
$K^{\rm nil}$ be the maximal extension of $K$ which for which the Galois group is a pro-$p$-group; here ${\rm nil}$ stands for nilpotent.  Then, as a pro-$2$-group, Serre proved in the early 1960s that 
\begin{equation}
\label{eq G2}
\Gal(\Q_2^{\rm nil}/\Q_2) = \langle x,y,z |x^2 y^4 x^{-1} y^{-1} x y \rangle.
\end{equation}
In general, Shafarevich proved in the 1940s that 
$\Gal(K^{\rm nil}/K)$ is free on $[K:\Q_p]+1$ generators
if $K$ does not contain $p^{\rm th}$ roots of unity.  Demushkin proved in the 1950s that it can always be presented with $[K:\Q_p]+2$ generators and one explicit relation, like \eqref{eq G2}.
}

For a finite group $G$, let $G/K$ be the set of 
Galois extensions $L^{\rm gal}/K$ in $\overline{\Q}_p/K$ with $\Gal(L^{\rm gal}/K) \cong G$.  These 
presentations let one compute 
cardinalities $|G/K|$.
Many completely explicit examples 
are given in \cite[\S 4]{Roe} for $K=\Q_p$.  For example, let $P_8$ be the Sylow
$2$-subgroup of $S_8$.  Then from \eqref{eq G2}  one eventually gets
$|P_8/\Q_2|=48$.  Switching language to our context of 
general extensions, let $nTj/K$ be the set of isomorphism classes of degree $n$ extensions
$L/K$ with $\Gal(L^{\rm gal}/K) = nTj$.   Then one can multiply by a constant
associated to $nTj$ to get the cardinality $|nTj/K|$.  For example,
over any field, each $P_8$ Galois extension comes from eight 
isomorphism classes of abstract $8T35$ extensions, and so $|8T35/\Q_2|=8 \cdot 48 = 384$.

The database lets one search by Galois group and see ramification
behavior not covered by the theorems.   Continuing the example of $8T35/\Q_2$, the database shows that they are 
distributed among eight families $2.1.8.C$ as summarized 
in Table~\ref{tab 8t35dist}.
\begin{table}[htb]
\[
\begin{array}{lrlll | r lll |  rrr}
C& M & s_1 & s_2 & s_3 &  m  & \multicolumn{3}{c}{\mbox{Hidden slopes}} &  P & \multicolumn{2}{c}{\mbox{Other $8Tj$}} \\
\hline
\ffamily{2.1.8.}{21a} & 32 & 1 & 1 & 2.75 & 16 & 2 & 2.5 & 2.5 &  2  & 8T38  \\
\ffamily{2.1.8.}{24d} & 16 & 1 & 2.5 & 2.75 & 16 & 1 & 2 & 2.5 & 1  \\
\hline
\ffamily{2.1.8.}{25b} & 32 & 1 & 2 & 3.25 & 16 & 1 & 2.5 & 3  & 3  & 8T21 & 8T31 \\ 
\ffamily{2.1.8.}{26b} & 32 & 1 & 2.5 & 3.25 & 16  &1 & 2 & 3 & 3 & 8T19 & 8T29  \\ 
\hline
\ffamily{2.1.8.}{27a} & 64 & 1 & 2.5 & 3.5 & 32 &  1& 2& 2.5 \mbox{ or } 3 & 7 & \multicolumn{2}{c}{\mbox{(five groups)}}\\ 
\hline
\ffamily{2.1.8.}{29a} & 64 & 2 & 2.5 & 3.75 & 32 &  1 & 3 & 3.25 & 3  & 8T28 & 8T30\\ 
\ffamily{2.1.8.}{30a} & 64 & 2 & 3 & 3.75 & 32 &  1 & 2.5 & 3.25 & 3 & 8T27 & 8T30  \\  
\hline
\ffamily{2.1.8.}{31a} & 128 & 2 & 3 & 4 & 32 &   1& 2.5 & 3.25  & 14 & \multicolumn{2}{c}{\mbox{(ten groups)}}  \\
\hline
&&&&&192&&&&
\end{array}
\]
\caption{\label{tab 8t35dist}The distribution of the 384 octic extensions of $\Q_2$ with associated Galois group $8T35$ into
eight families $2.1.8.C$}  
\end{table}
The table gives some indication of the nature of each of the families, including 
its total mass $M$ and the mass $m$ coming from $8T35$ fields.  
Each
$8T35$ field has mass $\frac{1}{2}$ and the sum $192$ of the $m$ is indeed $384/2$. 
The part of the family consisting of $8T35$ fields is always easy to 
describe.  For example, it is the locus where the coefficient $b_4=0$ in the
first-listed family $21a$.  When the number $P$ of packets is at most $3$, 
the remaining mass is evenly split among the remaining possible groups $8Tj$. 

In general, the organization of $p$-adic fields into families provides 
a framework for further 
investigation into ramification.  There are many resolvent maps  
$n'Tj'/K \rightarrow nTj/K$ coming from Galois theory 
over arbitrary fields.   The coordinates 
$d_0$, $a_\sigma$, $b_\sigma$, and $c_\sigma$ can be used to 
describe these maps in a  concise and uniform way.   Returning to the example, 
consider an $8T35$ extension
$L=K[x]/f(x)$.  If $f(x)$ is generically chosen then the degree
twenty-eight resolvent corresponding to the subgroup 
$S_6 \times S_2$ of $S_8$ factors into irreducibles as 
$f_4(x) f_8(x) f_{16}(x)$. The extension $L'=K[x]/f_8(x)$ is one of $L$'s seven different
siblings, the others then being easily obtainable either
by a degree 35 resolvent construction or by 
certain quadratic twists
\cite[Figure~3.1]{octic}.
The horizontal lines in Table~\ref{tab 8t35dist} indicate that for $K=\Q_2$ 
one has family interchanges $21a \leftrightarrow 24d$, $25b \leftrightarrow 26b$ and $29a \leftrightarrow 30a$ under this operation $L \leftrightarrow L'$.  Note that the set of 
six wild slopes is preserved in each of the three family interchanges as it must be.
The $8T35$ parts of the families $27a$ and $31a$ are closed under operation $L \leftrightarrow L'$.

\subsection{Making cyclic cases explicit via jump sets}
\label{subsec jumpsets} The natural goal of identifying the filtered group
$\Gal(\overline{\Q}_p/K)$ was reached nearly a century ago at the 
much simpler level of understanding the filtered abelianization 
$\Gal(K^{\rm ab}/K)$.  Namely local class field theory 
identifies $\Gal(K^{\rm ab}/K)$  
with the profinite completion of the multiplicative group $K^\times$ 
with the inertia group being sent to the unit group $U = \cO^\times$.  
For $j$ a positive integer and $s \in (j,j+1]$, the group
 $\Gal(K^{\rm ab}/K)^s$ 
is sent to the $j$-unit group $U_j = 1 + \Pi^j$.  

However this theoretically ideal solution does not immediately
answer some very basic concrete questions.  One such question
is, {\em what is the set $C_w(K)$ 
of Herbrand invariants $[s_1,\dots,s_w]$ coming from 
cyclic extensions of $K$ of degree $p^w$?}  In other words,
for what families $[s_1,\dots,s_w]/K$ is $(G,HSC) = (C_{p^w},[ \; ])$ one of the packets.
 One certainly needs the $s_k$ to form a strictly increasing
sequence of positive integers.   But to go beyond this statement, one 
needs to understand the filtered group $U_1$.  

Suppose the $p$-primary torsion in $U_1$ has order $p^v$.  
Then $U_1$ is isomorphic to $(\Z/p^v) \times \Z_p^{n_K}$.  
Thus the free vs.\ one-relator distinction from the beginning
of \S\ref{subsec givenGalois} is visible at this abelian level.   In the free case
$v=0$, the set $C_w(K)$ depends only on $e_K$ and is given below. 
 In the one-relator case $v\geq 1$,
the situation is much more complicated and $C_w(K)$ depends on an 
invariant $j_K$ extracted from the abelianization of the relation 
given in \eqref{eq mikirel} below.   The database tabulates
$j_K$, with instances having been given in Table~\ref{tab sixteentab}.  

The description of $C_w(K)$ involves combinatorial notions, as follows.  
For a prime $p$ and and a positive integer $e$, define
$\rho_{p,e}: \Z_{\geq 1} \rightarrow \Z_{\geq 1}$
 by $\rho_{p,e}(i) = \min(p i,i+e)$.   Table~\ref{tab intcolumns} draws 
 $\rho_{p,e}$ in three cases by organizing $\Z_{\geq 1}$ into
 $e$ columns.  Always $\rho_{p,e}(i)$ is the number immediately above 
 $i$ in its column.  Let $T_{p,e}$ be the set of non-images of $\rho_{p,e}$, 
 thus the $e$ numbers underlined and in bold at the bottom of columns.  
 A {\em jump set} of length $w$ is a sequence $[s_1,\dots,s_w]$ 
 satisfying the initial condition $s_1 \in T_{p,e}$ and
an inductive condition.  The inductive condition requires
that for $k \geq 2$ one has $s_k \geq \rho_{p,e}(s_{k-1})$, with
$s_k \in T_{p,e}$ if strict inequality holds.   Let $J_w(p,e)$ 
be the set of jump sets of length $w$.  If $p-1$ divides $e$ there
is also a notion of {\em extended jump set}.  Here $T_{p,e}$ is simply 
replaced by $T_{p,e}^* = T_{p,e} \cup \{pe/(p-1)\}$, the
extra point being indicated by bold italic in Table~\ref{tab intcolumns}.    

The set $J_w(p,e)$ of jump sets and the set $J^*_w(p,e)$ of
extended jump sets can be understood as the set of ways of climbing 
a ``rock wall.''  There are choices at the very beginning, but as soon as one reaches $e_K/(p-1)$, marked
by a light band in the table, the rest of the climbing path is forced to be vertical.  
The diagram splits into a part beneath $e_Kp/(p-1)$, marked by a dark band,
and the part on or above $e_Kp/(p-1)$.  We call the lower part 
{\em geometric}, because it agrees with the case $(p,\infty)$ corresponding
to $\F_q((t))$, and the upper part {\em arithmetic}.  
Table~\ref{tab intcolumns} also gives some sample counts, for $w=0$, $1$, $2$, \dots .
The count for the last printed $w$ also holds for all subsequent $w$.   

\begin{table}[htb]
\[
{\renewcommand{\arraystretch}{.2}
\renewcommand{\arraycolsep}{1pt}
{\renewcommand{\arraystretch}{.32}
\begin{array}{|cccccccc|}
 %\multicolumn{8}{c}{\;} \\
&&\;\;&&&\;\;&&20\\
&&&&&&19&\\
&18&&&&&&\\
&&&&17&&&\\
&&&16&&&&\\
15&&&&&&&\\
&&&&&&&14\\
&&&&&&13&\\
&12&&&&&&\\
&&&&11&&&\\
&&&10&&&&\\
\hline
\cellcolor[gray]{0.80} \boldsymbol{\mathit{9}}& \cellcolor[gray]{0.80} & \cellcolor[gray]{0.80}& \cellcolor[gray]{0.80}& \cellcolor[gray]{0.80}&\cellcolor[gray]{0.80}&\cellcolor[gray]{0.80}& \cellcolor[gray]{0.80}\\
\hline
&&&&&&&\underline{\mathbf{8}}\\
&&&&&&\underline{\mathbf{7}}&\\
&6&&&&&&\\
&&&&\underline{\mathbf{5}}&&&\\
&&&\underline{\mathbf{4}}&&&&\\
\hline
\cellcolor[gray]{0.90} 3&\cellcolor[gray]{0.90}&\cellcolor[gray]{0.90}&\cellcolor[gray]{0.90}&\cellcolor[gray]{0.90}&\cellcolor[gray]{0.90}&\cellcolor[gray]{0.90}&\cellcolor[gray]{0.90}\\
\hline
&\underline{\mathbf{2}}&&&&&&\\
\underline{\mathbf{1}}&&&&&&&\\
\hline
\multicolumn{8}{|c|}{(p,e) = (3,6)} \\
\hline
  \multicolumn{8}{c}{\;} \\
\multicolumn{8}{c}{|J_w(3,6)|=6,12} \\
  \multicolumn{8}{c}{\;} \\
\multicolumn{8}{c}{|J^*_w(3,6)|=7,15} \\
\end{array}
}
\;\;\;\;\;\;
\begin{array}{|cccccccc|}
&20&&&&&&\\
&&&&&19&&\\
&&18&&&&&\\
&&&&17&&&\\
\hline
\cellcolor[gray]{0.80} \boldsymbol{\mathit{16}}&\cellcolor[gray]{0.80}&\cellcolor[gray]{0.80}&\cellcolor[gray]{0.80}&\cellcolor[gray]{0.80}&\cellcolor[gray]{0.80}&\cellcolor[gray]{0.80}& \cellcolor[gray]{0.80}\\
\hline
&&&&&&&\underline{\mathbf{15}}\\
&&&14&&&&\\
&&&&&&\underline{\mathbf{13}}&\\
&12&&&&&&\\
&&&&&\underline{\mathbf{11}}&&\\
&&10&&&&&\\
&&&&\underline{\mathbf{9}}&&&\\
\hline
\cellcolor[gray]{0.90} 8&\cellcolor[gray]{0.90}&\cellcolor[gray]{0.90}&\cellcolor[gray]{0.90}&\cellcolor[gray]{0.90}&\cellcolor[gray]{0.90}&\cellcolor[gray]{0.90}&\cellcolor[gray]{0.90}\\
\hline
&&&\underline{\mathbf{7}}&&&&\\
 &6&&&&&&\\
  &&\underline{\mathbf{5}}&&&&&\\
 4&&&&&&&\\
  &\underline{\mathbf{3}}&&&&&&\\
 2&&&&&&&\\
\underline{\mathbf{1}}&&&&&&&\\
 \hline
 \multicolumn{8}{|c|}{(p,e) = (2,8)} \\
\hline
 \multicolumn{8}{c}{\;} \\
\multicolumn{8}{c}{|J_w(2,8)|=8,24,42,53} \\
 \multicolumn{8}{c}{\;} \\
\multicolumn{8}{c}{|J^*_w(2,8)|=9,23,53,69} \\
\end{array}
\;\;\;\;\;\;
\begin{array}{|ccccccccc|}
&&&&&&&20&\\
&&&&&&19&&\\
18&&&&&&&&\\
&&&&&17&&&\\
&&&&16&&&&\\
&15&&&&&&&\\
&&&14&&&&&\\
\hline
&&&&&&&&\underline{\mathbf{13}}\\
&&12&&&&&&\\
&&&&&&&\underline{\mathbf{11}}&\\
&&&&&&\underline{\mathbf{10}}&&\\
9&&&&&&&&\\
&&&&&\underline{\mathbf{8}}&&&\\
&&&&\underline{\mathbf{7}}&&&&\\
&6&&&&&&&\\
&&&\underline{\mathbf{5}}&&&&&\\
\hline
&&\underline{\mathbf{4}}&&&&&&\\
3&&&&&&&&\\
&\underline{\mathbf{2}}&&&&&&&\\
\underline{\mathbf{1}}&&&&&&&&\\
 \hline
 \multicolumn{9}{|c|}{(p,e) = (3,9)} \\
 \hline
 \multicolumn{9}{c}{\;} \\
 \multicolumn{9}{c}{|J_w(3,9)|=9,22,26} \\
\multicolumn{9}{c}{\;} \\
\multicolumn{9}{c}{\;} \\
\multicolumn{9}{c}{\;} \\
\multicolumn{9}{c}{\;} \\
\multicolumn{9}{c}{\;} \\
\multicolumn{9}{c}{\;} \\
\end{array}
}
\]
\caption{\label{tab intcolumns} Arrangements of the positive integers into $e$ columns 
for understanding jump sets. Bands are actually just convenient thickenings of 
the lines at heights $e_K/(p-1)$ and $e_Kp/(p-1)$, and they are drawn as lines
when these heights are not integral.} 
\end{table}

The desired identification for $v=0$ is $C_w(K)=J_w(p,e_K)$.  
As a trivial example, for $p>2$ one has the familiar 
$C_w(\Q_p) = J_w(p,1) = \{[1,\dots,w]\}$.   A partial answer for 
$v \geq 1$ is that $C_w(K) \subseteq J^*_w(p,e_K)$.   
As a simple example, $J^*_w(\Q_2) = \{[1,2,\dots,w],[2,3,\dots,w+1]\}$.
While $C_1(\Q_2) = \{[1],[2]\}$ is all of $J^*_1(\Q_2)$, otherwise
$C_w(\Q_2) = \{[2,3,\dots,w+1]\}$.  

In the $v \geq 1$ case where $K$ contains a primitive $p^{\rm th}$ root of unity $\zeta_p$, 
the next step towards a complete
answer goes as follows.  Building on Hasse \cite[Ch. 15]{hasse-zt}, Miki \cite[Lem 17]{miki}
showed that one can write
\begin{equation}
\label{eq mikirel}
\zeta_p=\alpha_0^{p^w}\alpha_1^{p^{w-1}}\dots
\alpha_{w-1}^p\alpha_w
\end{equation}
with each $\alpha_i$ satisfying one of the following
conditions:
\begin{itemize}
\item $v_K(\alpha_i)<pe_K/(p-1)$ and
$p\nmid v_K(\alpha_i)$,
\item $v_K(\alpha_i)=pe_K/(p-1)$ and $\alpha_i$ is not a
$p$th power,
\item $\alpha_i=1$.
\end{itemize}
Pagano \cite[\S 1.2]{pagano} defined the extended jump set
associated to $K$ by setting $j_i=v_K(\alpha_i-1)$ for
all $0\le i\le w$ such that $\alpha_i\not=1$; these
values are independent of the choice of
$\alpha_0,\alpha_1,\dots,\alpha_w$.  Undefined values in the
sequence $j = [j_0,j_1,\dots,j_w]$ are filled in using
the recursion $j_{i+1}=\rho_{p,e_K}(j_i)$.  
The remaining steps are to use $u_K$ to identify $C_w(K)$ 
as a subset of $J^*_w(p,e_K)$.
This is difficult to make
explicit in general, but Pagano works out many special 
cases.  

Pagano \cite[Thm.\,1.11]{pagano} also gave a formula for $j_K$ in terms of 
a defining Eisenstein polynomial for certain cases 
where $p \neq 2$ and all slopes of $K$ are
less than one.
The tabulation of invariants in 
the database using the factorization of $\zeta_p$ suggests that 
$j_K$ might always be directly expressible in 
terms of Eisenstein polynomials.  For the case $[1,\dots,1]/\Q_2$
recall from \S\ref{subsec unisloped} that the generic polynomial
is written 
\[2((1+2c_{2^w})+v_1x^{2^w-2^{w-1}}+\dots
+v_{w-1}x^{2^w-2^1}+v_wx^{2^w-2^0})+x^{2^w},\]  
with $v_w=1$. For $0\le k\le w$ set
$V_k = \sum_{i=1}^k v_i \in \F_2$.  Then we expect
\[
j_k = \left\{
\begin{array}{ll}
2^{k+1}-1 & \mbox{if $V_k=0$}, \\
2j_{k-1} & \mbox{if $V_k=1$ and $k<w$}, \\
2^{w+1}+2 c_{2^w}-2 & \mbox{if $V_k=1$ and $k=w$}.
\end{array}
\right.
\]
We also observe the simple statement that 
$j_K = [j_i]$ with $j_i = \epsilon (p^i-p+1)/(p-1)$ 
when all slopes are greater than $1$.  For the case
$(\epsilon,p)=(1,2)$ this formula gives $j_i = 2^i-1$. 
In particular, if $w=4$ then the jump set is $[1,3,7,15,31]$.

\subsection{Nonabelian quotients with known filtrations}
\label{subsec known}
The complications resolved by jump sets in the previous
subsection only became serious for $w \geq 2$.  Applying the easy special case of $w=1$ to tame extensions $K'$ of a fixed field $K$ gives a large class of extensions of $K$ which are 
nonabelian and wild, but still very well understood.  This 
situation has been studied from the point of view of primitive extensions of
$K$ by Del Corso, Dvornicich, and Monge in \cite{DDM}. 

To briefly summarize \cite{DDM} with some more specificity about slopes, 
let $K^{\rm prim} \subset \overline{\Q}_p$ 
be the composita of all primitive extensions of $K$.     It 
contains the maximal tame extension $K^{\rm tame}$ of
$K$.   What makes the group $\Gal(K^{\rm prim}/K)$  tractable is that its wild inertia subgroup $\Gal(K^{\rm prim}/K^{\rm tame})$ 
has exponent $p$.   The wild slopes are
exactly the positive rational numbers 
less than $e_K p/(p-1)$ with numerator and denominator prime to $p$, 
and then $e_K p/(p-1)$ itself.  The former occur with
infinite multiplicity but the latter occurs just with multiplicity 
$1$. The quotient $\Gal(K^{\rm tame}/K)$ is the closure 
of a subgroup $\langle \tau, \sigma|\sigma \tau \sigma^{-1} = \tau^q \rangle$,
where $\tau$ generates tame inertia and $\sigma$ is a Frobenius 
element.  The full group is a semidirect product 
$\Gal(K^{\rm prim}/K^{\rm tame}) \rtimes \Gal(K^{\rm tame}/K)$ with
a known action.

Consider now families $I=[s,\dots,s]_\epsilon^f/K$ with just 
one visible wild slope, such as $[\frac{3}{2}]_5/\Q_3$,
$[2]_7/\Q_2$, and $[1,\dots,1]/\Q_p$ of \S\ref{subsec easylinear}, \S\ref{subsec hardlinear} and
\S\ref{subsec unisloped}, respectively.  The compositum
$K^I \subset K^{\rm prim}$ of 
the Galois closures of the fields $L/K$ in $I/K$ is governed by the group theory just summarized.
So the packets $(G,HSC)$ that can occur
can be group-theoretically calculated.  

As an example of how group
theory can explain otherwise mysterious patterns, consider the
orbits of multiplication by $p$ on $\Z/(p^\rho-1)\Z$, with representatives
taken in $[\frac{p^\rho-1}{p-1},\frac{p(p^\rho-1)}{p-1})$.    In the case of $p^\rho = 2^3$,
two orbits are $\{8,9,11\}$, and $\{10,12,13\}$.  
Corresponding to dropping to the bottom of a column in a $(p,p^\rho-1)$ table
like Table~\ref{tab intcolumns}, remove all factors of $p$.  The
orbits in the example then become $\{1,9,11\}$, and $\{3,5,13\}$. 
The two triples contain the numbers in the previously mysterious borders of 
Table~\ref{tab fourteenpackets}, where one has to divide by $7$ to get
the hidden slopes for $[2]_7/\Q_2$.  The general recipe involves dividing by $p^\rho-1$ at the end.
In the example of $[\frac{3}{2}]_5/\Q_3$ from \S\ref{subsec easylinear}, $p^r = 3^4$ and
the orbit $\{56, 88, 104, 72\}$ becomes, via the drop $72 = 2^3 3^2 \rightarrow 8$, the hidden slopes $\frac{7}{10}$, $\frac{11}{10}$, $\frac{13}{10}$,
and $\frac{1}{10}$ mysteriously appearing there.  In general, 
the relatively elementary nature of $\Gal(K^{\rm prim}/K)$ makes
us hopeful that the explicit description of packets in the families of
 \S\ref{subsec easylinear}, 
 \S\ref{subsec hardlinear}, 
 and
\S\ref{subsec unisloped} will ultimately be specializations of the same
uniform description for any $[s,\dots,s]_\epsilon^f$. 

As a second example of a nonabelian quotient with a known filtration, 
consider the group $\Gal(K^{\rm nil}/K)$
defined just before \eqref{eq G2} and let $\Gal(K^{{\rm nil},p}/K)$ be
its maximal quotient of nilpotency class $p-1$ and exponent $p$.  
  In the strongest result of its type, Abrashkin \cite{Ab17} has identified the filtration on $\Gal(K^{{\rm nil},p}/K)$ under the assumption that $K$ contains a primitive $p$th root of unity.  For $p \geq 3$, Abrashkin's result, when it applies, goes well past the local class field theory of the previous subsection. However it leaves a lot of cases uncovered, as for example the Sylow $p$-subgroup of $S_{p^2}$ already has nilpotency class $p$.  
Thus for most of the families in the database, there is not yet a theoretical description of packets, even in principle.

\subsection{Comparing families $I/K$ for varying $K$} 
\label{subsec comparing}
For $s$ a real number and $K$ a finite degree subfield of $\overline{\Q}_p$, let 
$K^s \subset \overline{\Q}_p$ be the union of subextensions $L/K$ with 
all relative slopes less than $s$.   Say that two such fields $K$ and $K'$ are
$j$-close if there is an isomorphism of finite rings $\cO_K/\Pi_K^j \rightarrow 
\cO_{K'}/\Pi_{K'}^j$.  As an example, take $K'$ to have the same residual 
cardinality $q$ as $K$, and with ramification index at least $e_K$.  Then 
$K$ and $K'$ are $e_K$-close, because both finite rings are isomorphic
to $\F_q[t]/t^{e_K}$. 
Deligne proved in \cite[Th.\,2.8]{De}
that if $K$ and $K'$ 
are $j$-close then $\Gal(K^s/K)$ and $\Gal({K'}^s/K')$ are isomorphic 
as filtered groups.  Thus many instances of the problem of describing
filtered Galois groups have the same answer, even if we do not as yet
know the answer. 

To get a particular isomorphism between the Galois groups, well-defined up to conjugation,
one needs to choose a particular isomorphism between the finite rings.   
Suppose the choices of uniformizers and 
residue representatives made in \S\ref{subsec setup} are compatible with this ring
isomorphism.  Then the extensions $L/K$ and $L'/K'$ given by
the same specialization of the generic polynomial correspond.  
While Deligne's proof is complicated, he presents its basic
idea in \cite[\S 1.3]{De} as being this correspondence of Eisenstein polynomials.  
In terms of the database, the relative parts of the family pages for $I/K$ and $I/K'$ look
extremely similar, for any $I$ with top slope less than $j$.   

As one of the simplest possible examples, consider the Herbrand invariant $I = [\frac{1}{3},\frac{1}{3}]$ over $2$-adic fields. 
 Its generic polynomial is $\pi + \pi a_1 x + x^4$.  Restricting 
attention to $K$ with residual cardinality $2$, there is just one 
field $L/K$ in the family $I/K$, the one obtained by setting $a_1=1$.   For $\Q_2$, the field 
$L_0 = \Q_2[x]/(2 + 2 x + x^4)$ has associated Galois group  $S_4$ with hidden slope content $[ \; ]_3^2$. 
So for any $K$, the relative Galois group is likewise $S_4$ and the relative hidden slope content is likewise 
$[\;]_3^2$.   When one passes to absolute invariants, one naturally gets strong dependence on $K$.  
For example, taking $K$ to be one of the six ramified extensions of $\Q_2$ yields
trivial behavior: $L$ is the compositum $K \otimes L_0$ with Herbrand invariant $[\frac{1}{3},\frac{1}{3},s]$, 
Galois group $S_4 \times C_2$, and hidden slope content $[\;]_3^2$.   Here
$s \in \{1,2\}$ comes from the Herbrand invariant $[s]$ of $K/\Q_2$.   The unique ramified cubic extension 
$K = \Q_2[x]/(2-x^3)$ of 
$\Q_2$ yields completely different behavior.  Here the unique field $L$ is \field{2.1.12.12a1.1} with Galois group 
$C_2^6.C_9.C_6$ and slope content $[\frac{1}{9},\frac{1}{9},\frac{1}{9},\frac{1}{9},\frac{1}{9},\frac{1}{9}]_9^6$.  

As a more complicated example, take $I = [\frac{1}{2},\frac{2}{3}]$ with generic polynomial
$\pi (1+ a_3 x^3  + a_5 x^5 + c_6 x^6) + x^9$.   The page for $\family{3.1.9.13b}  =  [\frac{1}{2},\frac{2}{3}]/\Q_3$ 
says that there are eight fields, as $c_6$ is required to be $0$ 
in the parameter list $(a_3,a_5,c_6)$ whenever $a_3 \neq a_5$.   
The parameters $(1,2,0)$ and $(2,1,0)$ give fields having the
same Galois closure, with Galois group 9T18 and Galois slope content
$[\frac{1}{2},\frac{1}{2},\frac{2}{3}]_2^2$.  The remaining six fields from $(a_4,a_4,c_6)$ all have
Galois group $9T20$ and Galois slope content $[\frac{1}{2},\frac{1}{2},\frac{2}{3}]_2^3$,
with the splitting field depending only on $a_4$.   Since all slopes are less than one, Deligne's comparison
theorem says that all the facts just summarized also hold over any $K$ with residual
cardinality $3$.   Going further, one can compute via resolvent constructions
that the compositum $\Q_3^I$ has $\Gal(\Q_3^I/\Q_3)$ 
with order $2^4 3^8$ and slope content 
$[\frac{1}{2},\frac{1}{2},\frac{1}{2},\frac{1}{2},\frac{5}{9},\frac{2}{3},\frac{2}{3}]_2^6$.   So any $\Gal(K^I/K)$ has
the same structure as a filtered group.  

To see \S\ref{subsec jumpsets} in the light of the comparison theorem, let
$C^s_w(K)$ be the subset of $C_w(K)$ consisting of Herbrand
invariants $[s_1,\dots,s_w]$ with $s_w < s$.  Deligne's 
theorem says that $C_w^{s} (K) = C_w^{s} (\F_p((t)))$ 
holds for $s=e_K$, but in the easy case $v=0$ the explicit descriptions
says it holds for the larger integer $s = pe_K/(p-1)$.
It fails for even larger integers $s$, as one has entered the arithmetic regime.   Similarly, consider $I = [s,\dots,s]_e^f$.  Then Deligne's theorem says the filtered group
$\Gal(K^I/K)$ is isomorphic with its geometric analog
$\Gal(\F_q((t))^I/\F_q((t)))$ if $s < e_K$.  In the case $v=0$, the 
explicit description of $\Gal(K^I/K)$ indicated in \S\ref{subsec known} 
says that the isomorphism holds for $s<pe_K/(p-1)$.

 \subsection{Canonical globalization}   
 \label{subsec globalization} Our  complete tabulation of degree sixteen extensions of 
$\Q_2$ shows that the number of fields with automorphism group of size $1$, $2$, $4$, $8$,
and $16$ is respectively 9080, 833736, 44752, 2292, and 251.  These newly-determined numbers 
sum to the previously known total 890111.  The complete tabulation also gives the corresponding counts within each family. For the most ramified family \family{2.1.16.79a} $ = [2,3,4,5]$, 
the numbers are  $0$, $63488$,  $968$,  $240$, and $32$, for a total of $67728$.
  
The next step in our general approach to populating the database is to determine the Galois group $G$ and the
 hidden slope content $HSC$ for each field.  This is an ongoing process: we have computed 
 almost all the Galois groups using Doris' programs \cite{Do}, 
 but identifying $HSC$ is harder.  
 We describe here a completed part 
 that is of particular interest in terms of applications to number fields.  

For $r$ an odd prime, let $\Q^r \subset \C$ be the union of all finite degree 
Galois extensions  of $\Q$ with Galois group having order a power of $2$ and ramification within $\{\infty,2,r\}$.   Then a special case \cite[Example~11.18]{KO} of theorems of Koch gives two related
very strong results on $\Q^r$ for $r \equiv 3,5$ modulo $8$.  Let $D_\infty = \{1,c\}$ 
where $c$ is complex conjugation.  Let $D_r$ be the corresponding
decomposition group at the prime $r$, namely $D_r = \langle \tau,\sigma| \sigma \tau \sigma^{-1} = \tau^r \rangle$
understood in the category of pro-$2$-groups, so that $D_r$ has a semidirect product structure
$\Z_2 \rtimes \Z_2$. The first result is that $\Gal(\Q^r/\Q)$ is the free product $D_\infty * D_r$ in the 
category of pro-$2$-groups.  The second result is that its $2$-decomposition subgroup $D_2$ is
the entire global Galois group $\Gal(\Q^r/\Q)$.  

For $r \equiv 3,5 \; (8)$ fixed, Koch's result says that an extension 
$L/\Q_2$ with $G = \Gal(L^{\rm gal}/\Q_2)$ a $2$-group
globalizes either $0$ or $1$ times to a number field with discriminant $\pm 2^a r^b$ 
and Galois group $G$.   Typically a given field has no globalizations, as the governing
local group \eqref{eq G2} has three generators and one relation, while the governing 
global group also has three generators, but now has two relations, $c^2 = 1$ and 
$\sigma \tau \sigma^{-1} = \tau^r$.   

For $w$ a positive integer, consider the subfield $\Q^{r,w}$ generated by 
the subfields of $\Q^r$ of degree dividing $2^w$.  Its $r$-decomposition 
group is a quotient $D_{r,w}$ of $D_r$ having $2^{2w}$ elements and depending
only on $r$ mod $2^w$.  So the abstract group $\Gal(\Q^{r,w}/\Q) = D_\infty * D_{r,w}$ 
depends only on $r$ mod $2^w$.  Our calculations 
here show that the globalizing $2$-adic fields for $3$ and $19$ agree 
and the globalizing fields for $13$ and $29$ also agree.  This 
extends an observation made in \cite[\S8]{jones-roberts-number} for octic $2$-adic fields
and we do not have a proof that it holds for all $r$ and $w$.  

Table~\ref{tab stats} illustrates the frequency of globalization.  It counts the 
fields that  globalize exactly for the extra prime $r$ being in the indicated subset $R$ of 
$\{3,5,11,13\}$.  There are $1131$ $2$-adic fields that globalize to fields with discriminant 
$\pm 2^a$, of which $274$ have $a=79$.  So the bottom right entries, 
corresponding to  $R=\{3,5,11,13\}$, include many more always-globalizing 
fields.

\begin{table}[htb]
\[
{
\renewcommand{\arraycolsep}{4pt}
\begin{array}{r|rrrr}
& \multicolumn{4}{c}{\mbox{All $677795$ fields}} \\
  & \emptyset & \{3\} & \{11\} & \{3,11\}  \\
\hline
\emptyset  &   505520    & 9952 & 9952 & 61158  \\
\{5\} & 15072 & 0      &  0      & 0          \\
\{13\} & 15072 & 0      &  0      & 0          \\
\{5,13\} & 58076 & 0      & 0       & 2993   \\
\end{array} 
\;\;\;\;\;
\begin{array}{r|rrrr}
& \multicolumn{4}{c}{\mbox{Fields in $[2,3,4,5]/\Q_2$}} \\
  & \emptyset & \{3\} & \{11\} & \{3,11\}  \\
\hline
\emptyset  & 50614     & 1888 & 1888 &   8734 \\
\{5\} & 256    & 0       &  0      & 0          \\
\{13\} & 256   & 0       &  0      & 0         \\
\{5,13\} & 3384 & 0       & 0       &  708  \\
\end{array}
}
\]
\caption{\label{tab stats} Statistics of globalization for degree sixteen extensions
$L/\Q_2$ with Galois group $\Gal(L^{\rm gal}/\Q_2)$ having $2$-power order.}
\end{table}

    The purely-local relevance of globalization is that it allows easier
mechanical computation of Galois groups and it facilitates the identification
of hidden slopes.   Both $G$ and $HSC$ are in the database for all the canonically globalizing
fields just discussed.   Generally speaking, the database is designed so that it can present partial results.  
 As of this writing, the database shows $156$ packets inside the family $[2,3,4,5]$.   The 
numbers appearing as hidden wild slopes so far are $2$, $3$, $3.5$, $4$, $4.25$, $4.5$, $4.75$, $5.125$, $5.25$, $5.375$, 
and $5.625$.  

   Let $\Q_2^{{\rm nil},w}$  be the subfield of $\Q_2^{\rm nil}$ generated by subfields of degree $2^w$, so that $\Gal(\Q_2^{{\rm nil},w}/\Q_2)$ is a finite quotient of the infinite group $\Gal(\Q_2^{\rm nil}/\Q_2)$ of \eqref{eq G2}. 
Write its order as $2^{j_w}$.  A group-theoretical calculation says that $(j_1,j_2,j_3,j_4) = (3,8,25,204)$.    Corresponding numbers for 
$\Gal(\Q^{r,w}/\Q)$ begin independently of $r$, being $(3,7,18)$.   For $w=4$, there is dependence on $r$ mod $8$, with $p=3$ and $5$ yielding
$97$ and $101$.  

We have long known the twenty-five slopes appearing in compositum $\Q_2^{{\rm nil},3}$ of all nilpotent octic extensions of $\Q_2$.  There are three $-1$'s coming from the unramified octic extension of $\Q_2$, which is famously known not to globalize.  The wild slopes are then
    \[
    \mbox{$1$, $1$, $1$, $1\frac{1}{2}$, $1 \frac{1}{2}$, $\; 2$, $2$, $2$, $2\frac{1}{2}$, $2\frac{1}{2}$, $2\frac{5}{8}$, $2\frac{3}{4}$, $\;3$, 3,  $3\frac{1}{4}$,   $3\frac{1}{4}$, $3\frac{3}{8}$, $3\frac{1}{2}$, $3\frac{1}{2}$, $3\frac{3}{4}$, $3\frac{3}{4}$,   $\; 4$.}
    \]
     The eighteen slopes that survive
    to the quotient group $\Gal(\Q^{r,3}/\Q)$ are given for $r=3,5$ in \cite[\S8]{jones-roberts-number}. 

        The basic reason for constructing a large database is that many concrete facts about ramification are not yet known and seem resistant to theoretical investigation.  
        To underscore that much is not known, we conclude by asking a very concrete question:
      {\em what are the $204$ slopes of $\Gal(\Q_2^{{\rm nil},4}/\Q_2)$?}
    Basic theory says that four of them are $-1$, and the rest are in $[1,5]$ with 
    all denominators being powers of $2$.
    Analyzing high-degree composita of the new degree sixteen fields in the database will give many of these $204$ numbers.  We expect that a complete answer to the question may be out of reach without further theoretical advances, but computational progress can be measured by the number of slopes found.

\section*{Acknowledgements}
\label{sec acknowledgements}

This project arose from a SQuaRE at the American Institute of Mathematics, and we are grateful for their hospitality.  Roe was supported by Simons Foundation grant 550033.  We appreciate the comments made by the anonymous referees, which have improved the quality of the paper.

\bibliographystyle{amsalpha2}
\bibliography{local}
\end{document}